\documentclass{article}
\usepackage{xcolor,rotating,epic,eepic}
\usepackage{amssymb}
\usepackage{dsfont}
\usepackage{textcomp}
\usepackage{amsmath}
\usepackage{pgf}
\usepackage{txfonts}
\usepackage[colorlinks=true]{hyperref}
\usepackage{graphicx}

\title{Hardy type derivations on generalised series fields}

\author{Salma Kuhlmann, Micka\"{e}l Matusinski.}
\newcommand{\adresse}{\par\bigskip\small \rm

Salma Kuhlmann:\par Universit\"{a}t Konstanz,\par Fachbereich Mathematik und Statistik\par
78457 Konstanz, Germany.\par
Email: salma.kuhlmann@uni-konstanz.de\par
Webpage: http://www.math.uni-konstanz.de/~kuhlmann/
\pars
Micka\"{e}l Matusinski:\par Universit\"{a}t Konstanz,\par Fachbereich Mathematik und Statistik\par 78457 Konstanz, Germany.\par
Email: mickael.matusin@gmail.com\par 
Webpage: http://sites.google.com/site/mickaelmatusinski/}

\newtheorem{thm}{Theorem}[section]    
\newtheorem{lem}[thm]{Lemma}          
\newtheorem{prop}[thm]{Proposition}    
\newtheorem{cor}[thm]{Corollary}          
\newtheorem{definition}[thm]{Definition}
\newenvironment{defn}{\begin{definition}\rm}{\end{definition}}
\newtheorem{rem}[thm]{Remark}
\newenvironment{remark}{\begin{rem}\rm}{\end{rem}}
\newtheorem{exemple}[thm]{Example}
\newenvironment{ex}{\begin{exemple}\rm}{\end{exemple}}
\newtheorem{notation}[thm]{Notation}
\newenvironment{nota}{\begin{notation}\rm}{\end{notation}}
\newcommand{\n}{\par\noindent}

\newcommand{\sn}{\par\smallskip\noindent}
\newcommand{\mn}{\par\medskip\noindent}

\newcommand{\pars}{\par\smallskip}
\newtheorem{prf}{\it {Proof.}}
\newenvironment{demonstration}{\begin{prf}\rm}{\hfill$\Box$\end{prf}}

\date{October 3, 2011}
\begin{document}
\maketitle
\begin{abstract}    
We consider the valued field $\mathds{K}:=\mathbb{R}((\Gamma))$ of
generalised series (with real coefficients and monomials in a
totally ordered multiplicative group $\Gamma\>$). We investigate how
to endow $\mathds{K}$ with a series derivation, that is a derivation
that satisfies some natural properties such as commuting with
infinite sums (strong linearity) and (an infinite version of)
Leibniz rule. We characterize when such a derivation is of Hardy
type, that is, when it behaves like differentiation of germs of real
valued functions in a Hardy field. We provide a necessary and
sufficent condition for a series derivation of Hardy type to be
surjective.
\end{abstract}


\section{Introduction}
In his seminal paper, I. Kaplansky established
\cite[Corollary, p. 318]{kap} that if a valued field $(K, v)$
has the same characteristic as its residue field, then $(K, v)$ is
analytically isomorphic to a subfield of a suitable \textsl{field of
generalised series} (for definitions and terminology, see Section
\ref{sect:defi}). Fields of generalised series are thus universal
domains for valued fields. In particular, real closed fields of
generalised series provide suitable domains for the study of real
algebra.

The work presented in the first part of this paper is motivated by
the following query: are fields of generalised series suitable
domains for the study of real {\it differential} algebra?  We
investigate in Section 3 how to endow a field of generalised series
(of characteristic 0) with a natural derivation $d$, namely a
 \textsl{series derivation} (see Definition
 \ref{defi:series-deriv}). In the finite rank case (see Definition \ref{defn: hahn group}), the construction of such derivations presents no difficulty, as is already noticed in \cite{matu:puiseux-diff_rg-fini}. For arbitrary rank, but under an additional assumption (*) on the monomial group, examples of such series derivations are given in \cite{vdd:asch:liouv-cl-H-fields}. See Remark \ref{rem:rg-fini} for details on these questions. In this paper, we treat the general case.

 Our investigation is based on the notion of \textsl{fundamental monomials}, which are in fact representatives of the various
 comparability classes of series (see Definition \ref{defi:compar}). We start with a map $d$ from these fundamental monomials to the field
 of series. The central object of investigation is to extend $d$ first to the group of monomials (via a strong version of Leibniz rule)
and then from the group of monomials to the field of series (via an
infinite version of linearity) so that we obtain a series
derivation. The main challenge in doing so is to keep control of the
resulting supports and coefficients of the
resulting series. The criterion that we obtain in Theorem
\ref{theo:deriv} is rather abstract, but we derive from it more explicit results (Corollaries \ref{coro:(H1'')}, \ref{coro:(H1')} and \ref{coro:deriv'}). These results are applied in Section \ref{sect:exemple} to obtain concrete examples.
 \mn
Hardy fields, i.e. \textsl{fields of germs of differentiable real functions
at infinity}, were introduced by G. H. Hardy (the field of Log-Exp functions:\cite{hardy},\cite{hardy:LEfunct}) as the
natural domain for the study of asymptotic analysis. They represent prime examples of valued differential fields. In a series of papers,
 M. Rosenlicht studied the valuation theoretic properties of these derivations.
  This algebraic approach has been resumed and enhanced by M. Aschenbrenner and L. van den Dries in the formal axiomatic setting of H-fields \cite{vdd:asch:liouv-cl-H-fields}. The motivation for the second part of our paper is to understand the  possible connection between generalised series fields and Hardy fields as differential valued fields. Continuing our investigations in Section 4, we study derivations (on fields of generalised series) that satisfy the valuative properties discovered by Rosenlicht for Hardy
fields, namely \textsl{Hardy type derivations} (Definition \ref{defi:hardy_deriv}). This terminology comes from the notion of\textsl{ Hardy type asymptotic couple} in \cite{rosenlicht:val_gr_diff_val2}. We obtain in Theorem \ref{theo:hardy-deriv} a necessary and sufficient condition on a series derivation  to be of Hardy type. In the
 last section, we derive a criterion, Corollary \ref{int}, for a series derivation of Hardy type to be surjective.

A derivation on the Logarithmic-Exponential series field \cite{dmm:LE-pow-series} and on the field of transseries \cite{vdh:transs_diff_alg} have been introduced and studied. Furthermore, it is explained in \cite{Schm01} how to lift a given (strongly linear and compatible with the logarithm) derivation on a field of transseries to its exponential extensions. In \cite{matu-kuhlm:hardy-deriv-EL-series}, we extend our investigations to study Hardy type derivations on Exponential-Logarithmic series fields. In a forthcoming paper, we plan to endow the field of surreal numbers \cite{conway_numb-games,gonshor_surreal} with a derivation of Hardy type.

\section{Preliminary definitions}\label{sect:defi}
In this section, we introduce the required terminology and
notations.  For ordered set theory, we refer to
\cite{rosen:lin-ord}. In particular, we will repeatedly use the
following easy corollary of Ramsey's theorem \cite[ex. 7.5 p.
112]{rosen:lin-ord}:
\begin{lem}\label{lemme:suite}
Let $\Gamma$ be a totally ordered set. Every sequence
$(\gamma_n)_{n\in\mathbb{N}}\subset \Gamma$ has an infinite
subsequence which is either constant, or strictly increasing, or
strictly decreasing.
\end{lem}
\begin{defn} \label{defn: hahn group} Let $(\Phi,\preccurlyeq)$ be a totally ordered set, that we
call the set of \textbf{fundamental monomials}. We consider the set $\mathrm{\textbf{H}}(\Phi)$
of formal products $\gamma$ of the form
\[\gamma=\prod_{\phi\in\Phi}\phi^{\gamma_\phi}\]
 where $\gamma_\phi \in\mathbb{R}$, and the \textbf{support} of
 $\gamma$
\[\textrm{supp}\ \gamma := \{\phi\in\Phi\ |\ \gamma_{\phi}\neq 0\}\]
is an anti-well-ordered subset of $\Phi$.
 \end{defn}
Multiplication of formal products is defined pointwise: for $\alpha, \beta \in \mathrm{\textbf{H}}(\Phi)$
\[\alpha \beta = \prod _{\phi \in \Phi} \phi ^{\>\alpha _{\phi}+
\beta_{\phi}}\>\] With this multiplication, $\mathrm{\textbf{H}}(\Phi)$ is an abelian group with identity $1$ (the product with empty support). We endow
$\mathrm{\textbf{H}}(\Phi)$ with the anti lexicographic ordering $\preccurlyeq$ which extends $\preccurlyeq$ of
 $\Phi$: \[\gamma \succ 1
\mbox{ if and only if } \gamma _{\phi} > 0\> \mbox{ for } \phi :=
\max (\textrm{supp}\ \gamma) \>.\] With this definition, we see that
$\phi \succ 1$ for all $\phi\in \Phi$. Thus, $\mathrm{\textbf{H}}(\Phi)$ is a totally ordered abelian group \cite{hahn:nichtarchim}, that we call
 the \textbf{Hahn group} over $\Phi$.
\sn 
 Hahn's embedding theorem \cite{hahn:nichtarchim} states that an ordered abelian group $\Gamma$ embedds into $\mathrm{\textbf{H}}(\Phi)$ where $\Phi$ is the order type of its isolated subgroups. 
\textit{From now on, we consider some totally ordered set $(\Phi,\preccurlyeq)$ and we fix \textbf{$\Gamma$ subgroup of $\mathrm{\textbf{H}}(\Phi)$ with $\Phi\subset \Gamma$}}.  The set $\Phi$ is also called the \textbf{rank} of $\Gamma$.
\sn 
For any $\gamma\neq 1$, we will refer to $\gamma_\phi$ as the \textbf{exponent} of $\phi$, and the additive group $(\mathbb{R}, +)$ as the
 \textbf{group of exponents} of the fundamental monomials.

 \begin{defn} \label {defn: leading fundamental}
 We define the \textbf{leading fundamental monomial} of $\gamma\in \Gamma\setminus\{1\}$ by $\mbox{ LF }(\gamma):=\max(\textrm{supp}\
 \gamma)\>.$ We set by convention $\mbox{ LF }(1):=1$. This map verifies the \textbf{ultrametric triangular inequality}:
\[\forall \alpha,\beta\in\Gamma,
\mbox{ LF }(\alpha\beta)\preccurlyeq\max\{\mbox{ LF }(\alpha),\mbox{
LF }(\beta)\}\] and
\[\mbox{ LF }(\alpha\beta)=\max\{\mbox{ LF }(\alpha),\mbox{ LF }(\beta)\} \mbox{ if } \mbox{ LF }(\alpha)\neq
\mbox{ LF }(\beta)\>.\] 
\sn 
We define the \textbf{leading exponent}
of $1\ne \gamma \in \Gamma$ to be the exponent of $\mbox{ LF
}(\gamma)$, and we denote it by $\mbox{ LE }(\gamma)$. For $\alpha
\in \Gamma$  we set $|\alpha| := \max(\alpha, 1/\alpha)$; and define
\textbf{sign}$(\alpha)$ accordingly.
\end{defn}

In the following lemma, we summarize further properties of the maps
$\mbox{ LF }$ and $\mbox{ LE }$, that we will use implicitly
throughout the paper.
\begin{lem}\label{lemme:LF}
\item  1)\ \ For any $\alpha, \beta\in\Gamma\>,
\alpha\prec\beta\Leftrightarrow
 \mbox{ LE }\left(\displaystyle\frac{\beta}{\alpha}\right) > 0 \>.$
 \item 2)\ \ For any $1\ne \alpha \in\Gamma\>$ we have $\mbox{ LF }(|\alpha|)
 = \mbox{ LF }(\alpha)\>$ and $\mbox{ LE }(|\alpha|)
 = |\mbox{ LE }(\alpha)|\>.$
\item 3)\ \ We define on $\Gamma$ a scalar exponentiation: $\gamma
^r=(\prod_{\phi\in\Phi}\phi^{\gamma_\phi})^r:=\prod_{\phi\in\Phi}\phi^{r\gamma_\phi}\>\>$
for $r\in \mathbb{R}$. We have $\mbox{ LF }(\gamma ^r\>) = \mbox{ LF
}(\gamma)\>\>$ and $\mbox{ LE }(\gamma ^r\>) = r\> \mbox{ LE
}(\gamma)\>\>$, for $r\ne 0$.
\item 4)\ \ For $\beta\ne 1\ne \alpha \in\Gamma\>$ we have
\[\mbox{ LF }(\alpha) = \mbox{ LF }(\beta) \Leftrightarrow  \mbox{ there exists } n\in\mathbb{N} \mbox{ such that }
 |\beta|\preccurlyeq |\alpha|^n \mbox{ and }
|\alpha|\preccurlyeq |\beta|^n\>.\]
\item 5)\ \ For $\alpha$, $\beta \in \Gamma$ with $1\prec |\alpha|
  \prec |\beta|$,  we have $\mbox{ LF }(\alpha) \preccurlyeq \mbox{ LF }(\beta)$.
  \item 6)\ \ For $\alpha$, $\beta \in \Gamma$ with sign($\alpha$) =
  sign($\beta$), we have
  $\mbox{ LF }(\alpha\beta)=\max\{\mbox{ LF }(\alpha),\mbox{ LF }(\beta)\}$.
  \item 7)\ \ For any $\alpha, \beta\in\Gamma\>$, if
$\>\>\mbox{ LF }(\displaystyle\frac{\beta}{\alpha}) \prec \mbox{ LF
}(\beta)\>\>$ then $\mbox{ LF }(\alpha) = \mbox{ LF }(\beta)$ and
$\mbox{ LE }(\alpha) = \mbox{ LE }(\beta)$. In particular
sign($\alpha$) = sign($\beta$).
  \end{lem}
\begin{defn} \label{defn: generalised series}
Throughout this paper,
$\mathds{K}=\mathbb{R}\left(\left(\Gamma\right)\right)$ will denote
the \textbf{generalised series field with coefficients in
$\mathbb{R}$, and monomials in $\Gamma$}.
 It is the set of maps
\begin{center}
$\begin{array}{lccl} a\ :\ &\Gamma &\rightarrow &\mathbb{R}\\
&\alpha&\mapsto&a_\alpha\\
\end{array}$
\end{center} such that $\textrm{Supp}\ a:=\{\alpha\in\Gamma\ |\ a_\alpha\neq 0\}$ is anti-well-ordered in $\Gamma$.
As usual, we write these maps
$a=\displaystyle\sum_{\alpha\in\mathrm{Supp}\ a}a_\alpha \alpha$,
and denote by $0$ the series with
 empty support.
\sn By \cite{hahn:nichtarchim}, this set provided with
component-wise sum and the following convolution product
\[(\displaystyle\sum_{\alpha\in\mathrm{Supp}\ a}a_\alpha
\alpha\>)\>\>(\displaystyle\sum_{\beta\in\mathrm{Supp}\ b}b_\beta
\beta\>)\> =
 \displaystyle\sum_{\gamma\in\Gamma}\>(\displaystyle\sum_{\alpha\beta=\gamma}
 a_\alpha b_\beta)\> \>\gamma\]
 is a field.
\begin{remark}
The results in this paper hold for the generalised series field with
coefficients in an arbitrary ordered field $\mathcal{C}$ containing
$\mathbb{R}$ (instead of $\mathbb{R}$).
\end{remark}
\sn For any series $0\ne a$, we define its \textbf{leading monomial}:
$\mbox{ LM }(a) := \max\left(\mathrm{Supp}\ a\right)\in\Gamma\>$
with the usual convention that $\mbox{ LM }(0):=0\prec\gamma\>,$ for
all $\gamma \in \Gamma$.
The map \[\mbox{ LM }\>:\> \mathds{K}\setminus \{0\} \rightarrow \Gamma\] is a (multiplicatively written) \textbf{field valuation}; it verifies the following properties : \\
\[\forall a,b\in\mathds{K}: \mbox{ LM }(a.b)=\mbox{ LM }(a).\mbox{ LM }(b) \]\mbox{ and the \textbf{ultrametric triangular
inequality}}
\[ \mbox{ LM }(a+b)\preccurlyeq \max\{\mbox{ LM }(a),\mbox{ LM }(b)\}\>\]  \[ \mbox{ with }  \mbox{ LM }(a+b) = \max\{\mbox{ LM }(a),\mbox{
LM }(b)\} \mbox{ if } \mbox{ LM }(a)\ne \mbox{ LM }(b)\>.\] We
define the \textbf{leading coefficient} of a series to be $\mbox{ LC
}(a):= a_{\mbox{ LM }(a)}\in \mathbb{R}$ (with the convention that
$\mbox{ LC }(0)=0$) and use it to define a \textbf{total ordering} on
$\mathds{K}$ as follows: \[\forall a\in\mathds{K},\ a\leq
0\Leftrightarrow \mbox{ LC }(a)\leq 0\>\] For nonzero
$a\in\mathds{K}$, the term $\mbox{ LC }(a)\mbox{ LM }(a)$ is called
the \textbf{leading term} of $a$, that we denote $\mbox{LT}(a)$.
\end{defn} We use the leading monomial to extend the
ordering $\preccurlyeq$ on $\Gamma$ to a \textbf{dominance relation} on $\mathds{K}$  in the sense of G.H. Hardy (see \cite[Introduction p. 3-4]{hardy} and the Definition \ref{defi:dom_rel} below),
also denoted by $\preccurlyeq$: \begin{center}
$\forall a,b\in\mathds{K},\ a\preccurlyeq b\Leftrightarrow \mbox{ LM
}(a)\preccurlyeq \mbox{ LM }(b)$
\end{center} 
\begin{defn}\label{defi:dom_rel}
Let $(K,\leq)$ be an ordered field. A \textbf{dominance relation} on $K$ is a binary relation $\preccurlyeq$ on $K$ such that for all
$a,b,c\in K$:\\
\indent(DR1) $0\prec 1$\\
\indent(DR2) $a\preccurlyeq a$\\
\indent(DR3) $a\preccurlyeq b$ and $b\preccurlyeq c\ \Rightarrow\ a\preccurlyeq c$\\
\indent(DR4) $a\preccurlyeq b$ or $b\preccurlyeq a$\\
\indent(DR5) $a\preccurlyeq b\ \Rightarrow\ ac\preccurlyeq bc$\\
\indent(DR6) $a\preccurlyeq c$ and $b\preccurlyeq c\ \Rightarrow\ a-b\preccurlyeq c$\\
\indent(DR7) $0\leq a\leq b\ \Rightarrow\ a\preccurlyeq b$
\end{defn}
Given $a$ and $b$ non zero elements of $\mathds{K}$, we define the corresponding equivalence relations thus :
\begin{center}
$\begin{array}{lclcl}
 a\ \textrm{and}\ b\ \textrm{are }\textbf{asymptotic}& \Leftrightarrow\ &a\asymp b &\Leftrightarrow&  \mbox{ LM
}(a)= \mbox{ LM }(b)\\
 a\ \textrm{and}\ b\ \textrm{are }\textbf{equivalent}& \Leftrightarrow\ &a\sim b &\Leftrightarrow&  \mbox{ LT }(a)= \mbox{ LT }(b)
\end{array}$
\end{center}
\begin{defn}
We denote by $\mathds{K}^{\preccurlyeq 1} := \mathbb{R}((\Gamma^{\preccurlyeq 1}))= \{a\in\mathds{K}\ |\
a\preccurlyeq 1\}$ the \textbf{valuation ring} of $\mathds{K}$.
Similarly, we denote by $\mathds{K}^{\prec 1} :=\mathbb{R}((\Gamma^{\prec 1}))= \{a\in\mathds{K}\ |\ a\prec 1\}$ the \textbf{maximal ideal} of
$\mathds{K}^{\preccurlyeq 1}$. We have $\mathds{K}^{\preccurlyeq 1}=
\mathbb{R}\oplus \mathds{K}^{\prec 1}$. Thus $\mathbb{R}$ is isomorphic to the \textbf{residue field} $\mathds{K}^{\preccurlyeq 1}/\mathds{K}^{\prec 1}$ of $\mathds{K}$. We denote by
$\mathds{K}^{\succ 1} := \mathbb{R}\left(\left(\Gamma^{\succ
1}\right)\right)$, the \textbf{subring of purely infinite series}.
This is an additive complement group of
 $\mathds{K}^{\preccurlyeq 1}$ in $\mathds{K}$, i.e. $\mathds{K}=\mathds{K}^{\preccurlyeq 1}\oplus \mathds{K}^{\succ 1}$.
\end{defn}
\mn Finally, we extend the notion of \textbf{leading fundamental
monomial} to $\mathds{K}\backslash\{0\}$:
\begin{center}
$\begin{array}{lccl}
\mbox{ LF }\ :&\mathds{K}\setminus\{0\}&\rightarrow&\Phi\cup\{1\}\\
&a&\mapsto&\mbox{ LF }(a):=\mbox{ LF }(\mbox{ LM }(a))
\end{array}\>.$
\end{center}
We use it to define the notion of \textbf{comparability} of two series:
\begin{defn}\label{defi:compar}
Let $a\succ 1,\ b\succ 1$ be two elements of $\mathds{K}$. $a$ and
$b$ are \textbf{comparable} if and only if $\mbox{ LF }(a)=\mbox{ LF
}(b)$.
\end{defn}
It is straightforward to verify that comparability is an equivalence
relation on $\mathds{K}$.

\section{Defining derivations on generalised series fields}\label{section:deriv}

%

The following definition as in \cite[Part
II, Ch.8, Sect.5]{fuchs:partial_ord} will be needed to deal with infinite sums of series.
\begin{defn}\label{defi:summ_fam}
Let $I$ be an infinite index set and $\mathcal{F}=(a_i)_{i\in I}$ be a family
of series in $\mathds{K}$. Then $\mathcal{F}$ is said to be
\textbf{summable} if the two following properties hold:
\begin{description}
    \item[(SF1)] the support of the family $\mathrm{Supp}\ \mathcal{F}:=\displaystyle\bigcup_{i\in I}\textrm{Supp}\ a_i$ is
   anti-well-ordered in $\Gamma$;
\item[(SF2)] for any $\alpha\in\mathrm{Supp}\ \mathcal{F}$, the set $S_\alpha:=\{i\in I\ |\ \alpha\in\textrm{Supp}\ a_i\}\subseteq I$ is
finite.
\end{description}
Write $a_i=\displaystyle\sum_{\alpha\in\Gamma}a_{i,\alpha}\alpha$,
and assume that $\mathcal{F}=(a_i)_{i\in I}$ is summable. Then
\begin{center}
$\displaystyle\sum_{i\in I} a_i:=\displaystyle\sum_{\alpha\in
\mathrm{Supp}\ \mathcal{F}}\displaystyle\>(\sum_{i\in
S_\alpha}a_{i,\alpha}\>)\alpha\in\mathds{K}$
\end{center}
is a well defined element of $\mathds{K}$ that we call the
\textbf{sum} of $\mathcal{F}$.
\end{defn}
We will use subsequently the following characterisation of summability.

\begin{lem}\label{lemma:summability}
Given an infinite index set $I$ and a family $\mathcal{F}=(a_i)_{i\in I}$
of series in $\mathds{K}$, then $\mathcal{F}$ is summable if and only if the two following properties hold:
\begin{description}
    \item[(i)] for any sequence of monomials
$ (\alpha_n)_{n\in\mathbb{N}}\subset\mathrm{Supp}\ \mathcal{F}$, 
$\exists N\in \mathbb{N}$ such that $\alpha_N\succcurlyeq\alpha_{N+1}$;
\item[(ii)] for any sequence of \emph{pairwise distinct} indices $(i_n)_{n\in\mathbb{N}}\subset I$, $\displaystyle\bigcap_{n\in\mathbb{N}}\mathrm{Supp}\ a_{i_n}=\emptyset$.
\end{description}
\end{lem}
\begin{demonstration}
 Given a family $\mathcal{F}=(a_i)_{i\in I}$, the statement (i) is classically equivalent to the ``anti-well-orderedness" of $\mathrm{Supp}\ \mathcal{F}$, which is (SF1) (see e.g. \cite[Proposition 3.3]{rosen:lin-ord}).

Now suppose that (SF2) holds. Consider a sequence of pairwise distinct indices $(i_n)_{n\in\mathbb{N}}\subset I$ and the corresponding sequence of series $(a_{i_n})_{n\in\mathbb{N}}$ in $\mathcal{F}$. If there was some monomial $\alpha\in \displaystyle\bigcap_{n\in\mathbb{N}}\mathrm{Supp}\ a_{i_n}$, the corresponding set $S_\alpha$ would contain all the $i_n$'s and therefore would be infinite. This contradicts (SF2).

Suppose that (ii) holds, and that (SF2) fails, i.e. that there exists a monomial $\alpha\in\mathrm{Supp}\ \mathcal{F}$ such that the set $S_\alpha$ is infinite. Then we can choose in this set an infinite sequence of pairwise distinct indices $(i_n)_{n\in\mathbb{N}}$. Therefore, $\alpha\in\mathrm{Supp}\ a_{i_n}$ for all $n$, which means that $\alpha \in \displaystyle\bigcap_{n\in\mathbb{N}}\mathrm{Supp}\ a_{i_n}$. This contradicts (ii).
\end{demonstration}

Given a family $\mathcal{F}=(a_i)_{i\in I}$
of series with $I$ infinite, we call \textbf{subfamily of}  $\mathcal{F}$ any family  $\mathcal{F}'=(a_i)_{i\in J}$ for some index set $J\subset I$. By the preceding lemma, we note that \textsl{the family $\mathcal{F}$ is summable if and only if every countably infinite subfamily (i.e. with $J$ infinite countable) is summable}.
\n We introduce in the following definition the precise notion of ``good" derivation for generalised series.
\begin{defn}\label{defi:series-deriv}
Given the generalised series field $\mathds{K}$, consider the following axioms:
\begin{description}
\item[(D0)]  $1'=0$;
\item[(D1) Strong Leibniz rule:]$\forall\alpha=\displaystyle\prod_{\phi\in\mathrm{supp}\
\alpha}\phi^{\alpha_\phi}\in\Gamma,\ \alpha '=
\alpha\displaystyle\sum_{\phi\in\mathrm{supp}\ \alpha}\alpha_\phi\displaystyle\frac {\phi '}{\phi}$;
    \item[(D2) Strong linearity:]$\forall a=
\displaystyle\sum_{\alpha\in\mathrm{Supp}\ a}a_\alpha
\alpha\in\mathds{K},\ a'= \displaystyle\sum_{\alpha\in\mathrm{Supp}\
a}a_\alpha \alpha '$.
\end{description}
A map $d:\Gamma\rightarrow\mathds{K}$, $\alpha\mapsto \alpha '$, verifying (D0) and (D1) is called a \textbf{series derivation on} $\Gamma$.
A map $d:\mathds{K}\rightarrow\mathds{K}$, $a\mapsto a'$, verifying these three axioms is called a \textbf{series derivation on }$\mathds{K}$.
\end{defn}
 \begin{remark}
\sn \ \ A series derivation is a derivation in the usual sense, i.\
e.:
\begin{enumerate}
\item $d$ is linear : $\forall a,b\in\mathds{K},\ \forall K,L\in\mathbb{R},\ (K.a+L.b)'=K.a'+L.b'$.
\item $d$ verifies the Leibniz rule : $\forall a,b\in\mathds{K},\ (ab)'=a'b+ab'$.
\end{enumerate}
\end{remark}
The problem arising from the preceding definition, which is the main purpose of this section, is to clarify when the axioms (D1) and (D2) make sense. More precisely, we want to characterise the existence of such series derivations by some specific properties of their restriction to fundamental monomials.
\begin{defn}\label{defi:series-deriv-extension}
Let
\begin{center}$\begin{array}{lccl}
d_\Phi\ :&\Phi&\rightarrow& \mathds{K}\backslash\{0\}\\
&\phi &\mapsto& \phi '
\end{array}$
\end{center}
be a map.\\
1) We say that $d_\Phi$ \textbf{extends to a series derivation on $\Gamma$} if the following property holds:
\begin{description}
\item[(SD1)] for any $\alpha\in\Gamma$,
the family $\left(\displaystyle\frac {\phi '}{\phi}\right)_{\phi\in \mathrm{supp}\ \alpha}$ is summable. \end{description} Then the \textbf{series
derivation} $d_\Gamma$ on $\Gamma$ (extending $d_\Phi$) is defined
to be the map \[d_\Gamma :\ \Gamma\rightarrow\mathds{K}\] obtained
through the axioms (D0) and (D1) (which clearly makes sense by (SD1)).
\mn 2) We say that a series derivation $d_\Gamma$ on $\Gamma$
\textbf{extends to a series derivation on $\mathds{K}$} if the
following property holds:
\begin{description}
\item[(SD2)] for any  $a\in\mathds{K}$, the family $(\alpha ')_{\alpha\in \mathrm{Supp}\ a}$ is summable.
\end{description}
Then the \textbf{series derivation $d$ on $\mathds{K}$ (extending
$d_\Gamma$)} is defined to be the map \[d\ :\
\mathds{K}\rightarrow\mathds{K}\] obtained through the
axiom (D2) (which clearly makes sense by (SD2)).
\end{defn}

\begin{remark}\label{rem:rg-fini}
\begin{enumerate}
\item As is already noticed in \cite[Definition 2.2]{matu:puiseux-diff_rg-fini}, when the fundamental chain $\Phi$ is finite, say $\Phi=\{\phi_1,\ldots, \phi_r\}$ for some $r\in\mathbb{N}^*$, then any map $d_\Phi\ :\Phi\rightarrow \mathds{K}\backslash\{0\}$ extends to a series derivation on $\Gamma$ and on $\mathds{K}$. Indeed:\begin{enumerate}
\item for any monomial $\alpha=\phi_1^{\alpha_1}\cdots\phi_r^{\alpha_r}\in\Gamma$, $\alpha '=
\alpha.\left(\alpha_1\displaystyle\frac {\phi_1 '}{\phi_1}+\cdots+\alpha_r\displaystyle\frac {\phi_r '}{\phi_r}\right)$ is well-defined ;
\item for any series $a=
\displaystyle\sum_{\alpha\in\mathrm{Supp}\ a}a_\alpha
\alpha\in\mathds{K}$, \begin{center}
$a'= \displaystyle\sum_{\alpha\in\mathrm{Supp}\
a}a_\alpha \alpha ' = \left(\displaystyle\sum_{\alpha\in\mathrm{Supp}\
a}a_\alpha \alpha_1. \alpha\right)\displaystyle\frac {\phi_1 '}{\phi_1}+\cdots+  \left(\displaystyle\sum_{\alpha\in\mathrm{Supp}\ a}a_\alpha \alpha_r. \alpha\right )\displaystyle\frac {\phi_r '}{\phi_r}$
\end{center} is well-defined.
\end{enumerate}
\item In \cite[Section 11]{vdd:asch:liouv-cl-H-fields}, the authors define a derivation $d$ on $\mathds{K}$ under the assumption that the monomial group $\Gamma$ satisfies a condition called (*) (i.e. admits a valuation basis; see \cite{kuhl:ord-exp}). In this case, $\Gamma \simeq \mathrm{\textbf{H}}_{\mathrm{fin}}(\Phi)$, the subgroup of $\mathrm{\textbf{H}}(\Phi)$ of monomials with finite support. So (SD1) is easily verified as in (a) above. We note that this derivation $d$ is what we call a \textbf{monomial derivation} (see Definition \ref{defi:monom-deriv}). In Section \ref{sect:exemple} we analyse how to obtain (SD2) in this monomial derivation case (see Proposition \ref{prop:monomial_case}).
\end{enumerate}
\end{remark}

 In the next Theorem \ref{theo:deriv}, we provide a necessary and sufficient condition on a
 map $d_\Phi :\Phi\rightarrow \mathds{K}$ so that properties (SD1) and (SD2) hold. In the sequel, we drop the subscripts $\Phi$ and $\Gamma$ of $d_\Phi$
 and $d_\Gamma$ to relax the notation. We isolate the following two crucial ''bad'' hypotheses:
\begin{description}
\item[{\bf(H1)}] there exists a strictly decreasing sequence $(\phi_n)_{n\in\mathbb{N}}\subset\Phi$ and an increasing sequence
 $(\tau^{(n)})_{n\in\mathbb{N}}\subset\Gamma$ such that for any $n$, $\tau^{(n)}\in\mathrm{Supp}\ \displaystyle\frac{\phi_n'}{\phi_n}$;\\
\item[{\bf(H2)}] there exist strictly increasing sequences $(\phi_n)_{n\in\mathbb{N}}\subset\Phi$ and
$(\tau^{(n)})_{n\in\mathbb{N}}\subset\Gamma$ such that for any $n$,
$\tau^{(n)}\in\mathrm{Supp}\ \displaystyle\frac{\phi_n'}{\phi_n}$
and $\mbox{ LF
}\left(\displaystyle\frac{\tau^{(n+1)}}{\tau^{(n)}}\right)\succcurlyeq\phi_{n+1}$.
\end{description}
\begin{thm}\label{theo:deriv}
 A map $d\ :\
\Phi\rightarrow \mathds{K}\backslash\{0\}$ extends to a series
derivation on $\mathbb{K}$ if and only both hypotheses $(H1)$ and
$(H2)$ fail.
\end{thm}
The proof of this theorem will be split into the proofs of Lemma \ref{lemme:H1} and Lemma \ref{lemme:H2}.

\begin{remark}\label{rem:H2'}
Let a series derivation $d$ on $\Gamma$ be given. We claim that the following condition (H2') is a positive version of (H2), i.e. a condition that will be necessary and sufficient for (SD2) to hold:
\begin{description}
\item[(H2')] for any strictly increasing sequences
$(\phi_n)_{n\in\mathbb{N}}\subset\Phi$ and
$(\tau^{(n)})_{n\in\mathbb{N}}\subset\Gamma$ such that for any $n$,
$\tau^{(n)}\in\mathrm{Supp}\ \displaystyle\frac{\phi_n'}{\phi_n}$,
the set $\mathcal{S}=\left\{n\in\mathbb{N}\ |\ \mbox{ LF
}\left(\displaystyle\frac{\tau^{(n+1)}}{\tau^{(n)}}\right) \succcurlyeq\phi_{n+1}\right\}$ is finite.
\end{description}
Indeed, the Hypothesis (H2') implies clearly that (H2) does not hold. 
Conversely, suppose that there exist strictly increasing sequences $(\phi_n)_{n\in\mathbb{N}}\subset\Phi$ and
$(\tau^{(n)})_{n\in\mathbb{N}}\subset\Gamma$ as in (H2'), for which
$\mathcal{S}$ is infinite. Denote $\mathcal{S}=\{n_i\ |\ i\in\mathbb{N}\}$ with $n_i<n_{i+1}$ for all $i$, and set $m_i:=n_i+1$, $i\in\mathbb{N}$. We notice that
$\displaystyle\frac{\tau^{(m_{i+1})}}{\tau^{(m_i)}}= \displaystyle\frac{\tau^{(n_{i+1}+1)}}{\tau^{(n_i+1)}}=
\displaystyle\frac{\tau^{(n_{i+1}+1)}}{\tau^{(n_{i+1})}}
\displaystyle\frac{\tau^{(n_{i+1})}}{\tau^{(n_{i+1}-1)}}\cdots
\displaystyle\frac{\tau^{(n_i+2)}}{\tau^{(n_i+1)}}$. Moreover we
have $\mbox{ LF
}\left(\displaystyle\frac{\tau^{(n_{i+1}+1)}}{\tau^{(n_{i+1})}}\right)\succcurlyeq
\phi_{n_{i+1}+1}$ and for any $n$ such that $n_i<n<n_{i+1}$, $\mbox{ LF
}\left(\displaystyle\frac{\tau^{(n+1)}}{\tau^{(n)}}\right)\prec\phi_{n+1}$.
So applying the ultrametric inequality for $\mbox{ LF }$ (see Definition \ref{defn: leading fundamental}), we have $\mbox{ LF
}\left(\displaystyle\frac{\tau^{(m_{i+1})}}{\tau^{(m_{i})}}\right)=
\mbox{ LF
}\left(\displaystyle\frac{\tau^{(n_{i+1}+1)}}{\tau^{(n_{i+1})}}\right) \succcurlyeq \phi_{n_{i+1}+1}=\phi_{m_{i+1}}$. Thus the increasing sequences $(\phi_{m_i})_{i\in\mathbb{N}}$ and $(\tau^{(m_i)})_{i\in\mathbb{N}}$ verify (H2).
\end{remark}

\mn To emphasise the role of each hypothesis, we divide the proof of the
Theorem \ref{theo:deriv} into the statement and the proof of the two following lemmas \ref{lemme:H1}, \ref{lemme:H2}.
\begin{lem}\label{lemme:H1}
A map $d\ :\Phi \rightarrow  \mathds{K}\backslash\{0\}$ extends to
 a series derivation on $\Gamma$ if and only if $(H1)$ fails.
\end{lem}
\begin{demonstration}
Suppose that $(H1)$ holds, i.e. there exist a strictly decreasing sequence
$(\phi_n)_{n\in\mathbb{N}}$ and an increasing one $(\tau^{(n)})_{n\in\mathbb{N}}$ such that for all $n$,
$\tau^{(n)}\in\mathrm{Supp}\ \displaystyle\frac{\phi _n '}{\phi_n}$. Applying Lemma \ref{lemme:suite} to the sequence $(\tau^{(n)})_{n\in\mathbb{N}}$, we have two possibilities. Either there is an increasing subsequence, which contradicts the point (i) of Lemma \ref{lemma:summability}. Or there is a constant one, which implies that
$\displaystyle\bigcap_{n\in\mathbb{N}}\mathrm{Supp }\displaystyle\frac{\phi _n '}{\phi_n}\neq\emptyset$, contradicting the point (ii) of Lemma \ref{lemma:summability}. Thus the family $\left(\displaystyle\frac{\phi _n '}{\phi_n}\right)_{n\in\mathbb{N}}$ is not summable.
\mn Conversely, suppose that $(SD1)$ does not hold. There exists an infinite
anti-well-ordered subset $E:=\mathrm{supp}\ \alpha\subset\Phi$ such that the family
$\left(\displaystyle\frac {\phi '}{\phi}\right)_{\phi\in E}$ fails
to be summable. By the Lemma \ref{lemma:summability}, there are two cases. Contradicting point (ii), there exists a sequence $(\phi_n)_{n\in\mathbb{N}}$ of pairwise distinct fundamental monomials so that there exists a monomial  $\tau\in\displaystyle\bigcap_{n\in\mathbb{N}}\mathrm{Supp}\ \displaystyle\frac{\phi _n '}{\phi_n}$. Then just define $\tau^{(n)}:=\tau$ for all $n$. Contradicting point (i), there exists a strictly increasing sequence of monomials $(\tau^{(n)})_{n\in\mathbb{N}}$ in $\displaystyle\bigcup_{\phi\in E}\mathrm{Supp}\ \displaystyle\frac{\phi '}{\phi}$. Subsequently, for any $n\in\mathbb{N}$, choose $\phi_n\in E$ so that $\tau^{(n)}\in\mathrm{Supp}\ \displaystyle\frac{\phi _n '}{\phi_n}$. Since it is a sequence from $E$ which is anti-well-ordered, $(\phi_n)_{n\in\mathbb{N}}$ cannot contain any strictly increasing subsequence. Moreover, we claim that, without loss of generality, the $\phi_n $'s may be assumed to be pairwise distinct. Indeed, since for any $
 \phi \in E$, $\mbox{ Supp }(\displaystyle\frac{\phi '}{\phi})$ is anti-well-ordered in $\Gamma$, the set
$\{\tau^{(n)} \>|\>n\in\mathbb{N}\} \cap \mbox{ Supp }(\displaystyle\frac{\phi '}{\phi})$ is finite. In other words, the
map
\begin{center}
$\begin{array}{lccl} &\{ \tau^{(n)} \>|\>
n\in\mathbb{N}\}&\rightarrow &\{ \phi_n \>|\>n\in\mathbb{N}\}
\\
&\tau^{(n)}&\mapsto&\phi_n
\end{array}\>$
\end{center}
has infinitely many finite fibres. Choosing a complete set of
representatives for the set of fibres, we may extract a subsequence of $(\tau^{(n)})_{n\in\mathbb{N}}$ (which is strictly increasing as is $(\tau^{(n)})_{n\in\mathbb{N}}$) and with pairwise distinct corresponding $\phi_n$'s. We continue to
denote such a subsequence by $(\tau^{(n)})_{n\in\mathbb{N}}$ below .

Now applying Lemma \ref{lemme:suite} to the sequence $(\phi_n)_{n\in\mathbb{N}}$, we obtain that it must contain a strictly decreasing subsequence. Such subsequence together with the corresponding $\tau^{(n)}$'s are the sequences complying the requirements of (H1).
\end{demonstration}

Now we introduce a new tool that will help us to derive from the preceding lemma  more concrete corollaries and several examples. Given an anti-well-ordered set $E$, we denote by $ot(E)$ its \textbf{order type} \cite{rosen:lin-ord}.
\begin{defn}\label{defi:isom}
\begin{itemize}
\item Consider $\mu,\nu\in\Phi$ such that $ot(\mathrm{Supp}\
\displaystyle\frac{\mu '}{\mu})\leq ot(\mathrm{Supp}\
\displaystyle\frac{\nu '}{\nu})$. There exists an isomorphism of ordered sets from $\mathrm{Supp}\
\displaystyle\frac{\mu '}{\mu}$ onto a final segment of $\mathrm{Supp}\
\displaystyle\frac{\nu '}{\nu}$. In the sequel, we shall denote this isomorphism by $I_{\mu,\nu}$, and its inverse isomorphism $I_{\mu,\nu}^{-1}$ by $I_{\nu,\mu}$. Note that $I_{\mu,\nu}(\textrm{LM}(\displaystyle\frac{\mu '}{\mu}))=\textrm{LM}(\displaystyle\frac{\nu '}{\nu})$.
\item Consider $\phi,\psi\in\Phi$.  We shall say that $I_{\mu,\nu}$ is a \textbf{left
shift} if $I_{\mu,\nu}(\gamma) \prec \gamma$ for any $\gamma$ in
the domain of $I_{\mu,\nu}$.
\item We can enumerate the elements of $\mathrm{Supp}\ \displaystyle\frac{\phi '}{\phi}$ in the decreasing
direction $\tau_0\succ\tau_1\succ\cdots\succ\tau_\lambda\succ\cdots$
where $\lambda$ is an ordinal number called the \textbf{position} of $\tau_\lambda$ in  $\mathrm{Supp}\ \displaystyle\frac{\phi '}{\phi}$. Thus, denoting \textbf{ON} the proper class of all ordinals \cite{rosen:lin-ord}, we define the set of functions $\{p_\phi,\phi\in\Phi\}$ by:
\begin{center}
$\forall\phi\in\Phi,\ \ p_\phi : \mathrm{Supp}\ \displaystyle\frac{\phi '}{\phi}\rightarrow \mathrm{\textbf{ON}}$
\end{center} which maps any element $\tau_\lambda\in\mathrm{Supp}\ \displaystyle\frac{\phi '}{\phi}$ to its position $\lambda$ in $\mathrm{Supp}\ \displaystyle\frac{\phi '}{\phi}$.\\
Note that, given any $\phi,\psi\in\Phi$ and any $\tau^{(\phi)},\tau^{(\psi)}$ in the domain of $I_{\phi,\psi}$, respectively $I_{\psi,\phi}$, we have $p_\phi(\tau^{(\phi)})=p_\psi(\tau^{(\psi)})$ if and only if $I_{\phi,\psi}(\tau^{(\phi)})=\tau^{(\psi)}$ (if and only if $I_{\psi,\phi}(\tau^{(\psi)})=\tau^{(\phi)}$).
\end{itemize}
\end{defn}

\begin{lem}\label{lemma:isom}
Suppose that (H1) holds (or equivalently that (SD1) does not hold). The corresponding strictly decreasing sequence $(\phi_n)_{n\in\mathbb{N}}$ from $\Phi$ and the increasing sequence $(\tau^{(n)})_{n\in\mathbb{N}}$ from $\Gamma$ with $\tau^{(n)}\in\mathrm{Supp}\ \displaystyle\frac{\phi_n'}{\phi_n}$, can be chosen so that for any $n\in\mathbb{N}$, $\tau^{(n)}$ is in the domain of $I_{\phi_n,\phi_{n+1}}$ and $\tau^{(n)}\preccurlyeq \tau^{(n+1)}\preccurlyeq I_{\phi_n,\phi_{n+1}}(\tau^{(n)})$. In particular, the sequence $(I_{\phi_n,\phi_{n+1}})_{n\in\mathbb{N}}$ consists in automorphisms that are not left shifts.
\end{lem}
\begin{demonstration}
Consider from (H1) a strictly decreasing sequence $(\phi_n)_{n\in\mathbb{N}}\subset\Phi$ and an increasing sequence
 $(\tau^{(n)})_{n\in\mathbb{N}}\subset\Gamma$ such that for any $n$, $\tau^{(n)}\in\mathrm{Supp}\ \displaystyle\frac{\phi_n'}{\phi_n}$. Consider $S:=\{p_{\phi_n}(\tau^{(n)}),\ n\in\mathbb{N}\}$ which is a subset of \textbf{ON} and for any $\lambda\in S$, consider $S_\lambda:=\{n\in\mathbb{N}\ |\ p_{\phi_n}(\tau^{(n)})=\lambda\}$ (see Definition \ref{defi:isom}).

Suppose that there exists $\lambda\in S$ such that $S_\lambda$ is infinite. So it contains a strictly increasing subsequence $(n_i)_{i\in\mathbb{N}}$ of natural numbers. Since the sequence $\tau^{(n)}$ is increasing by (H1), for any $i\in\mathbb{N}$, we have $\tau^{(n_i)}\preccurlyeq\tau^{(n_{i+1})}= I_{\phi_{n_i},\phi_{n_{i+1}}}(\tau^{(n_i)})$. The sequences $(\phi_{n_i})_{i\in\mathbb{N}}$ and $(\tau^{(n_i)})_{i\in\mathbb{N}}$ have the required properties. Suppose now that for any $\lambda\in S$, the set $S_\lambda$ is finite. This implies that $S$ is infinite. For any $m\in\mathbb{N}$, denote $S^{(m)}:=\{p_{\phi_n}(\tau^{(n)}),\ n>m\}$ and $S^{(m)}_\lambda:=\{n\in\mathbb{N}\ |\ n>m \textrm{ and } p_{\phi_n}(\tau^{(n)})=\lambda\}$.
We shall define by induction a strictly increasing sequence $(\lambda_i)_{i\in\mathbb{N}}$ from $S$, together with the desired sequence $(\tau^{(n_j)})_{j\in\mathbb{N}}$. Set $\lambda_0:=\min S$. Then denote $S_{\lambda_0}=\{n_0,\ldots,n_{j_0}\}$ with $n_{k+1}>n_k$ for any $k$. Consider the corresponding monomials $\tau^{(n_0)} \preccurlyeq\tau^{(n_1)}\preccurlyeq\cdots\preccurlyeq\tau^{(n_{j_0})}$. Since for any $k$, $p_{\phi_{n_k}}(\tau^{(n_k)})=\lambda_0$, we have $\tau^{(n_{k})}\preccurlyeq\tau^{(n_{k+1})}= I_{\phi_{n_k},\phi_{n_{k+1}}}(\tau^{(n_{k})})$ as desired.

Now suppose that we have built a finite sequence $\tau^{(n_0)}\preccurlyeq\tau^{(n_1)}\preccurlyeq\cdots \preccurlyeq\tau^{(n_{j_0})} \preccurlyeq\cdots\preccurlyeq\tau^{(n_{j_i})}$ together with an ordinal $\lambda_i$ for some $i\geq 0$, with the desired properties. Then, set $\lambda_{i+1}:=\min S^{(n_{j_i})}$, which implies that $\lambda_{i+1}>\lambda_i$ (all the indices $n$ corresponding to lower ordinals $\lambda_l$ are lower than $n_{j_i}$). Now consider the set $S^{(n_{j_i})}_{\lambda_{i+1}}$ which is non empty by definition of $\lambda_{i+1}$. Then we denote it $S^{(n_{j_i})}_{\lambda_{i+1}}=\{n_{j_i+1},\ldots,n_{j_{i+1}}\}$ with $n_{j_i+k+1}>n_{j_i+k}$ for any $k$. Then the corresponding monomials are such that $\tau^{(n_{j_i+k})}\preccurlyeq \tau^{(n_{j_i+k+1})}=I_{\phi_{n_{j_i+k}},\phi_{n_{j_i+k+1}}}(\tau^{(n_{j_i+k})})$ for any $k$. Moreover, since $n_{j_i}<n_{j_i+1}$ and $\lambda_i<\lambda_{i+1}$, we have $\tau^{(n_{j_i})}\preccurlyeq \tau^{(n_{j_i+1})}\prec I_{\phi_{n_{j_i}},\phi_{n_{j_i+1}}}(\tau^{(n_{j_i})})$ as desired.
\end{demonstration}

We deduce from the preceding lemma a more explicit sufficient
condition (but not necessary: see Example \ref{ex:(H'2)}) such that
(SD1) holds.
\begin{cor}\label{coro:(H1')} A map $d\ :\ \Phi\rightarrow\mathds{K}\backslash\{0\}$ extends to a series derivation on $\Gamma$ if the following property holds :\begin{description}
\item[(H1')]the set $E_1=\left\{\phi\in\Phi\ |\ \exists\ \psi\succ\phi,\
I_{\psi,\phi}\textrm{ is not a left shift}\right\}$ is well ordered
in $\Phi$.
\end{description}
\end{cor}
\begin{demonstration} For any strictly decreasing sequence $S=(\phi_n)_{n\in\mathbb{N}}$, since $E_1\subset\Phi$ is well-ordered,
 $E_1\cap S$ is finite. So all but finitely many couples $(\phi_n,\phi_{n+1})$ are such that $I_{\phi_n,\phi_{n+1}}$ is a left shift. It implies that we can not obtain a sequence $(\tau^{(n)})_{n\in\mathbb{N}}$ as in (H1).
\end{demonstration}
To visualize (H1'), we illustrate in the following Figure 1, the supports
$\mathrm{Supp}\ \displaystyle\frac{\phi '}{\phi}$ for some
$\phi\in\Phi$. The ordered sets $\Phi$ and $\Gamma$ are represented
as linear orderings.

\begin{figure}
\begin{pgfpicture}{-5.4781cm}{-5.4781cm}{5.5219cm}{5.5219cm}%
\pgfsetxvec{\pgfxy(0.8264,0)}
\pgfsetyvec{\pgfxy(0,0.8264)}
\pgfsetroundjoin%
\pgfsetstrokecolor{black}
\pgfsetlinewidth{0.2pt} 
\pgfputat{\pgfxy(0.5,-4.4)}{\pgftext{\color{black}\large $\Gamma$}}\pgfstroke
\pgfsetlinewidth{0.8pt} 
\pgfsetfillcolor{black}\pgfcircle[fillstroke]{\pgfxy(0.5,-3.5)}{.025cm}
\pgfcircle[fillstroke]{\pgfxy(0.5,-3)}{.025cm}
\pgfcircle[fillstroke]{\pgfxy(0.5,-2.5)}{.025cm}
\pgfcircle[fillstroke]{\pgfxy(0.5,3.5)}{.025cm}
\pgfcircle[fillstroke]{\pgfxy(0.5,4)}{.025cm}
\pgfcircle[fillstroke]{\pgfxy(0.5,4.5)}{.025cm}
\pgfcircle[fillstroke]{\pgfxy(-3.5,-1)}{.05cm}
\pgfcircle[fillstroke]{\pgfxy(-3.4,-1)}{.05cm}
\pgfcircle[fillstroke]{\pgfxy(-3.2,-1)}{.05cm}
\pgfcircle[fillstroke]{\pgfxy(-2.7,-1)}{.05cm}
\pgfcircle[fillstroke]{\pgfxy(-2,-1)}{.05cm}
\pgfcircle[fillstroke]{\pgfxy(-1.1,-1)}{.05cm}
\pgfcircle[fillstroke]{\pgfxy(0.1,-1)}{.05cm}
\pgfcircle[fillstroke]{\pgfxy(0.2,-1)}{.05cm}
\pgfcircle[fillstroke]{\pgfxy(0.4,-1)}{.05cm}
\pgfcircle[fillstroke]{\pgfxy(0.9,-1)}{.05cm}
\pgfcircle[fillstroke]{\pgfxy(1.6,-1)}{.05cm}
\pgfcircle[fillstroke]{\pgfxy(2.5,-1)}{.05cm}
\pgfcircle[fillstroke]{\pgfxy(-3.2,1)}{.05cm}
\pgfcircle[fillstroke]{\pgfxy(-3.1,1)}{.05cm}
\pgfcircle[fillstroke]{\pgfxy(-2.9,1)}{.05cm}
\pgfcircle[fillstroke]{\pgfxy(-2.4,1)}{.05cm}
\pgfcircle[fillstroke]{\pgfxy(-1.7,1)}{.05cm}
\pgfcircle[fillstroke]{\pgfxy(-0.8,1)}{.05cm}
\pgfcircle[fillstroke]{\pgfxy(0.4,1)}{.05cm}
\pgfcircle[fillstroke]{\pgfxy(0.5,1)}{.05cm}
\pgfcircle[fillstroke]{\pgfxy(0.7,1)}{.05cm}
\pgfcircle[fillstroke]{\pgfxy(1.2,1)}{.05cm}
\pgfcircle[fillstroke]{\pgfxy(1.9,1)}{.05cm}
\pgfcircle[fillstroke]{\pgfxy(2.8,1)}{.05cm}
\pgfcircle[fillstroke]{\pgfxy(-3,-0)}{.05cm}
\pgfcircle[fillstroke]{\pgfxy(-2.9,-0)}{.05cm}
\pgfcircle[fillstroke]{\pgfxy(-2.7,-0)}{.05cm}
\pgfcircle[fillstroke]{\pgfxy(-1.5,-0)}{.05cm}
\pgfcircle[fillstroke]{\pgfxy(0.8,0)}{.05cm}
\pgfcircle[fillstroke]{\pgfxy(0.9,0)}{.05cm}
\pgfcircle[fillstroke]{\pgfxy(1,0)}{.05cm}
\pgfcircle[fillstroke]{\pgfxy(1.2,0)}{.05cm}
\pgfcircle[fillstroke]{\pgfxy(1.5,0)}{.05cm}
\pgfcircle[fillstroke]{\pgfxy(1.9,0)}{.05cm}
\pgfcircle[fillstroke]{\pgfxy(2.4,0)}{.05cm}
\pgfcircle[fillstroke]{\pgfxy(3,0)}{.05cm}
\pgfcircle[fillstroke]{\pgfxy(-2.9,2)}{.05cm}
\pgfcircle[fillstroke]{\pgfxy(-2.8,2)}{.05cm}
\pgfcircle[fillstroke]{\pgfxy(-2.6,2)}{.05cm}
\pgfcircle[fillstroke]{\pgfxy(-2.1,2)}{.05cm}
\pgfcircle[fillstroke]{\pgfxy(-1.4,2)}{.05cm}
\pgfcircle[fillstroke]{\pgfxy(-0.5,2)}{.05cm}
\pgfcircle[fillstroke]{\pgfxy(0.7,2)}{.05cm}
\pgfcircle[fillstroke]{\pgfxy(0.8,2)}{.05cm}
\pgfcircle[fillstroke]{\pgfxy(1,2)}{.05cm}
\pgfcircle[fillstroke]{\pgfxy(1.5,2)}{.05cm}
\pgfcircle[fillstroke]{\pgfxy(2.2,2)}{.05cm}
\pgfcircle[fillstroke]{\pgfxy(3.1,2)}{.05cm}
\pgfsetlinewidth{0.2pt} 
\pgfputat{\pgfxy(-4.1,1)}{\pgftext[right]{\color{black}\small $\psi\in(\Phi\backslash E_1)$}}\pgfstroke
\pgfsetarrows{-to}
\pgfmoveto{\pgfxy(3.07,1.9)}\pgflineto{\pgfxy(2.83,1.1)}
\pgfstroke
\pgfsetarrows{-}
\pgfputat{\pgfxy(3.05,1.5)}{\pgftext[left]{\color{black}\small $I_{\phi,\psi}$}}\pgfstroke
\pgfsetlinewidth{0.8pt} 
\pgfputat{\pgfxy(-4.8,4)}{\pgftext{\color{black}\large $\Phi$}}\pgfstroke
\pgfsetlinewidth{0.2pt} 
\pgfputat{\pgfxy(-4.1,2)}{\pgftext[right]{\color{black}\small $\phi\in(\Phi\backslash E_1)$}}\pgfstroke
\pgfputat{\pgfxy(-4.1,-0)}{\pgftext[right]{\color{black}\small $\tilde{\phi}\in E_1$}}\pgfstroke
\pgfsetarrows{-to}
\pgfmoveto{\pgfxy(-4,-4)}\pgflineto{\pgfxy(6,-4)}
\pgfstroke
\pgfsetarrows{-}
\pgfsetarrows{-to}
\pgfmoveto{\pgfxy(-4,-4)}\pgflineto{\pgfxy(-4,6)}
\pgfstroke
\pgfsetarrows{-}
\pgfputat{\pgfxy(0.5,-5.5)}{\pgftext{\color{black}\small Figure 1. Illustration of $(H1')$}}\pgfstroke
\end{pgfpicture}%
\end{figure}

Under an additional hypothesis, we deduce from Lemma \ref{lemma:isom}
a necessary and sufficient condition for a map $d$ on $\Phi$ to extend to a series derivation on $\Gamma$:
\begin{cor}\label{coro:(H1'')}
 Let a map $d\ :\ \Phi\rightarrow\mathds{K}\backslash\{0\}$ be given. We suppose that there exists $N\in\mathbb{N}$ such that, for any $\phi\in\Phi$, $\mathrm{Card}\left(\mathrm{Supp}\ \displaystyle\frac{\phi '}{\phi}\right)\leq N$. Then $d$ extends to a series derivation on $\Gamma$ if and only if the following property holds :\\
\textbf{(H1'')} for any strictly decreasing sequence
$(\phi_n)_{n\in\mathbb{N}}\subset\Phi$, there exists a pair of
integers $m<n$ such that $I_{\phi_{m},\phi_{n}}$ is a left shift.
\end{cor}
\begin{demonstration}
 Suppose that (SD1) does
not hold. Equivalently, by (H1),
there exist a strictly decreasing sequence $(\phi_n)_{n\in\mathbb{N}}$ and an
increasing one
$(\tau^{(n)})_{n\in\mathbb{N}}$ with $\tau^{(n)}\in\mathrm{Supp}\
\displaystyle\frac{\phi_n'}{\phi_n}$ for any $n$.
We set $k_n:=p_{\phi_n}(\tau^{(n)})\in\{1,\ldots,N\}$, $n\in\mathbb{N}$ (see Definition \ref{defi:isom}).
Applying Lemma \ref{lemme:suite} to the sequence
$(k_n)_{n\in\mathbb{N}}$, there exists an infinite constant subsequence
$(k_{n_i}=k)_{i\in\mathbb{N}}$. Hence, for any $i<j$,
$\tau^{(n_i)}\preccurlyeq\tau^{(n_{j})}=I_{\phi_{n_i},\phi_{n_{j}}}(\tau^{(n_i)})$ (see the final remark in Definition \ref{defi:isom}).
The sequence $(\phi_{n_i})_{i\in\mathbb{N}}$ is such that the
corresponding isomorphisms $I_{\phi_{n_i},\phi_{n_{j}}}$ for any
$i<j$ fail to be left shifts.
\mn Conversely, suppose that there
exists a decreasing sequence $(\phi_n)_{n\in\mathbb{N}}$ for
which the $I_{\phi_m,\phi_{n}}$'s, $m<n$, are not left shifts. That
is, given $m$, for any $n>m$, there exists $\tau^{(m)}\in\mathrm{Supp}\
\displaystyle\frac{\phi_m '}{\phi_m}$ such that
$\tau^{(m)}\preccurlyeq I_{\phi_m,\phi_{n}}(\tau^{(m)})$. Thus for any $n$, we set $l^{(m)}_{n}:=p_{\phi_m}(\tau^{(m)})\in\{1,\ldots,N\}$. By Lemma \ref{lemme:suite}, there exists a constant subsequence $(l^{(m)}_{n_i}=l^{(m)})_{i\in\mathbb{N}}$, that is we have
$\tau^{(m)}\preccurlyeq
I_{\phi_m,\phi_{n_i}}(\tau^{(m)})$ for any $i\geq 0$. Now, consider the sequence
$(l^{(m)})_{m\in\mathbb{N}}$. Again by Lemma \ref{lemme:suite} and since for any $m$ $l^{(m)}\in\{1,\ldots,N\}$, there exists a constant subsequence, say
$(l^{(m_j)}=l)_{j\in\mathbb{N}}$ for some $l\in\{1,\ldots,N\}$. Hence for any $j\in\mathbb{N}$, $\tau^{(m_j)}\preccurlyeq
I_{\phi_{m_j},\phi_{m_{j +1}}}(\tau^{(m_j)})=\tau^{(m_{j+1})}$. The sequence $(\tau^{(m_j)})_{j\in\mathbb{N}}$ verifies (H1), which means that (SD1) does not hold for the family $\{\phi_{m_j}\ |\ j\in\mathbb{N}\}$.
\end{demonstration}

\begin{ex}\label{ex:(H'2)}
In Corollary \ref{coro:(H1'')}, the assumption that the cardinalities of the sets $\mathrm{Supp}\ \displaystyle\frac{\phi '}{\phi}$, $\phi\in\Phi$, are uniformely bounded is necessary to apply the criterion (H1''). Indeed, if we drop this assumption, (SD1) may still hold even if (H1'') fails, as illustrated by the following Figure 2. The dashed lines indicate changes of comparability classes (for
instance, take $\tau_{0,k}=\phi_1\prec
I_{\phi_0,\phi_k}(\tau_{0,k})=\phi_0^{1/k}$ for any
$k\in\mathds{N}^*$). The lines connect $\tau_{k,l}$ and
$I_{\phi_k,\phi_l}(\tau_{k,l})$ for which the isomorphism
$I_{\phi_k,\phi_l}$ fails to be a left shift.\\
We observe that, even if there is an infinite decreasing sequence $(\phi_n)_{n\in\mathbb{N}}$
for which the $I_{\phi_n,\phi_{n+1}}$'s are not left shifts, (SD1)
holds for the anti-well-ordered subset $\{\phi_n,\ n\in\mathbb{N}\}$ of $\Phi$. Indeed, by construction, the set $\displaystyle\bigcup_{n\in\mathbb{N}}\mathrm{Supp}\ \displaystyle\frac{\phi_n '}{\phi_n}$ is anti-well-ordered and $\displaystyle\bigcap_{n\in\mathbb{N}}\mathrm{Supp}\ \displaystyle\frac{\phi_n '}{\phi_n}=\emptyset$ (see Lemma \ref{lemma:summability}).

\begin{figure}
\begin{pgfpicture}{-5.4781cm}{-5.4781cm}{5.5219cm}{5.5219cm}%
\pgfsetxvec{\pgfxy(0.8264,0)}
\pgfsetyvec{\pgfxy(0,0.8264)}
\pgfsetroundjoin%
\pgfsetstrokecolor{black}
\pgfsetlinewidth{0.2pt} 
\pgfputat{\pgfxy(0.5,-4.4)}{\pgftext{\color{black}\large $\Gamma$}}\pgfstroke
\pgfsetarrows{-to}
\pgfmoveto{\pgfxy(-4,-4)}\pgflineto{\pgfxy(6,-4)}
\pgfstroke
\pgfsetarrows{-}
\pgfsetarrows{-to}
\pgfmoveto{\pgfxy(-4,-4)}\pgflineto{\pgfxy(-4,6)}
\pgfstroke
\pgfsetarrows{-}
\pgfsetfillcolor{black}\pgfcircle[fillstroke]{\pgfxy(5,3)}{.05cm}
\pgfcircle[fillstroke]{\pgfxy(4.8,2)}{.05cm}
\pgfcircle[fillstroke]{\pgfxy(4.6,1)}{.05cm}
\pgfcircle[fillstroke]{\pgfxy(4.4,0)}{.05cm}
\pgfcircle[fillstroke]{\pgfxy(4.2,-1)}{.05cm}
\pgfcircle[fillstroke]{\pgfxy(3.2,4)}{.05cm}
\pgfcircle[fillstroke]{\pgfxy(2.8,2)}{.05cm}
\pgfcircle[fillstroke]{\pgfxy(2.6,1)}{.05cm}
\pgfcircle[fillstroke]{\pgfxy(2.4,0)}{.05cm}
\pgfcircle[fillstroke]{\pgfxy(2.2,-1)}{.05cm}
\pgfcircle[fillstroke]{\pgfxy(1,3)}{.05cm}
\pgfcircle[fillstroke]{\pgfxy(0.6,1)}{.05cm}
\pgfcircle[fillstroke]{\pgfxy(0.4,0)}{.05cm}
\pgfcircle[fillstroke]{\pgfxy(0.2,-1)}{.05cm}
\pgfcircle[fillstroke]{\pgfxy(-1.2,2)}{.05cm}
\pgfcircle[fillstroke]{\pgfxy(-1.6,-0)}{.05cm}
\pgfcircle[fillstroke]{\pgfxy(-1.8,-1)}{.05cm}
\pgfcircle[fillstroke]{\pgfxy(0.65,-3)}{.05cm}
\pgfcircle[fillstroke]{\pgfxy(0.65,-2.5)}{.05cm}
\pgfcircle[fillstroke]{\pgfxy(0.65,-2)}{.05cm}
\pgfputat{\pgfxy(-4.3,4)}{\pgftext{\color{black}\small $\phi_0$}}\pgfstroke
\pgfputat{\pgfxy(-4.3,3)}{\pgftext{\color{black}\small $\phi_1$}}\pgfstroke
\pgfputat{\pgfxy(-4.3,2)}{\pgftext{\color{black}\small $\phi_2$}}\pgfstroke
\pgfcircle[fillstroke]{\pgfxy(-4.3,0.5)}{.025cm}
\pgfcircle[fillstroke]{\pgfxy(-4.3,-0)}{.025cm}
\pgfcircle[fillstroke]{\pgfxy(-4.3,-0.5)}{.025cm}
\pgfsetdash{{5pt}{3pt}}{0pt}
\pgfmoveto{\pgfxy(3.4,4.5)}\pgflineto{\pgfxy(3.4,-1.5)}
\pgfstroke
\pgfmoveto{\pgfxy(1.2,4.5)}\pgflineto{\pgfxy(1.2,-1.5)}
\pgfstroke
\pgfmoveto{\pgfxy(-1,4.5)}\pgflineto{\pgfxy(-1,-1.5)}
\pgfstroke
\pgfsetdash{}{0pt}
\pgfmoveto{\pgfxy(3.2,4)}\pgflineto{\pgfxy(5,3)}
\pgflineto{\pgfxy(3.2,4)}\pgflineto{\pgfxy(4.8,2)}
\pgflineto{\pgfxy(3.2,4)}\pgflineto{\pgfxy(4.6,1)}
\pgflineto{\pgfxy(3.2,4)}\pgflineto{\pgfxy(4.4,0)}
\pgflineto{\pgfxy(3.2,4)}\pgflineto{\pgfxy(4.2,-1)}
\pgfstroke
\pgfmoveto{\pgfxy(1,3)}\pgflineto{\pgfxy(2.8,2)}
\pgflineto{\pgfxy(1,3)}\pgflineto{\pgfxy(2.6,1)}
\pgflineto{\pgfxy(1,3)}\pgflineto{\pgfxy(2.4,0)}
\pgflineto{\pgfxy(1,3)}\pgflineto{\pgfxy(2.2,-1)}
\pgfstroke
\pgfmoveto{\pgfxy(-1.2,2)}\pgflineto{\pgfxy(0.6,1)}
\pgflineto{\pgfxy(-1.2,2)}\pgflineto{\pgfxy(0.4,0)}
\pgflineto{\pgfxy(-1.2,2)}\pgflineto{\pgfxy(0.2,-1)}
\pgfstroke
\pgfputat{\pgfxy(0.5,-5.2)}{\pgftext{\color{black}\small Figure 2. Counter-example when Card(Supp  $\frac{\phi  '}{\phi}$)}}\pgfstroke
\pgfputat{\pgfxy(0.5,-5.7)}{\pgftext{\color{black}\small is not uniformely bounded.}}\pgfstroke
\pgfputat{\pgfxy(-5.5,1)}{\pgftext{\color{black}\large $\Phi$}}\pgfstroke
\end{pgfpicture}%
\end{figure}
\end{ex}
\mn Now we prove the second lemma that completes the proof of
Theorem \ref{theo:deriv}.
\begin{lem}\label{lemme:H2}
Let $d$ a series derivation on $\Gamma$ be given. Then $d$ extends
to a series derivation on $\mathds{K}$ if and only if (H2) fails.
\end{lem}
\begin{demonstration}
First, we suppose that (H2) holds. For any $n\in\mathbb{N}$, set $\displaystyle\frac{\tau^{(n+1)}}{\tau^{(n)}}=\psi_{n+1}^{\eta_{n+1}} \gamma^{(n+1)}$ where
$\psi_{n+1}=\mbox{ LF
}\left(\displaystyle\frac{\tau^{(n+1)}}{\tau^{(n)}}\right)$, $\eta_{n+1}=\mbox{ LE
}\left(\displaystyle\frac{\tau^{(n+1)}}{\tau^{(n)}}\right)$ and $\gamma^{(n+1)}\in\Gamma$.
Then $\psi_{n+1}\succcurlyeq\phi_{n+1}$, $\eta_{n+1}>0$ (the
sequence $(\tau^{(n)})_n$ is strictly increasing) and $\mbox{ LF
}\left(\gamma^{(n+1)}\right)\prec\psi_{n+1}$. Consider now the sequence $(\alpha^{(n)})_{n\in\mathbb{N}}$ where $\alpha^{(0)}=\phi_0^{-\epsilon_0}$ for some $\epsilon_0>0$,
$\alpha^{(n+1)}=\phi_{n+1}^{-\epsilon_{n+1}}$ for some
$\epsilon_{n+1}>0$ if $\psi_{n+1}\succ\phi_{n+1}$, and $\alpha^{(n+1)}=(\phi_{n+1}^{\eta_{n+1}}\gamma^{(n+1)})^{-1}= \displaystyle\frac{\tau^{(n)}}{\tau^{(n+1)}}$ if $\psi_{n+1}=\phi_{n+1}$. This
sequence is decreasing since the sequence
$(\phi_n)_{n\in\mathbb{N}}$ is increasing. Moreover, setting  $\beta^{(n)}:=\alpha^{(n)}\tau^{(n)}$, we
have $\beta^{(n)}\in\mathrm{Supp}\ (\alpha^{(n)})'$ for any $n$ (see (D1): $\phi_n\in\mathrm{Supp }\alpha^{(n)}$ and $\tau^{(n)}\in\mathrm{Supp }\displaystyle\frac{\phi_n'}{\phi_n}$). Then it is routine to prove that  $\displaystyle\frac{\beta^{(n+1)}}{\beta^{(n)}}= \displaystyle\frac{\alpha^{(n+1)}\tau^{(n+1)}}{\alpha^{(n)}\tau^{(n)}}\succ 1$, meaning the sequence $(\beta^{(n)})_{n\in\mathbb{N}}$ is strictly increasing. It implies that the family $((\alpha^{(n)})')_{n\in\mathbb{N}}$ is not summable, witnessing that (SD2) does not hold.
\n 
Conversely, suppose that (SD2) does
not hold, i.e. there exists an anti-well-ordered set of monomials
$E\subset\Gamma$ such that the family $(\alpha ')_{\alpha\in E}$ is not summable. By Lemma \ref{lemma:summability}, there are two cases.
Either the set
 $\displaystyle\bigcup_{\alpha\in E}\mathrm{Supp}\ \alpha '$ contains a strictly increasing sequence $(\beta^{(n)})_{n\in\mathbb{N}}$, or there exists a subsequence $(\alpha^{(n)})_{n\in\mathbb{N}}$ of pairwise distinct elements of $E$ such that $\beta\in\displaystyle\bigcap_{n\in\mathbb{N}}\mathrm{Supp}\ (\alpha^{(n)})'$ for some $\beta\in\Gamma$. In the latter case, we denote (as in the former) by $\beta^{(n)}=\beta$ some copy of $\beta$ in $\mathrm{Supp}\ (\alpha^{(n)})'$: the sequence $(\beta^{(n)})_{n\in\mathbb{N}}$ is constant. 

In the former case, set $(\alpha^{(n)})_{n\in\mathbb{N}}$ a corresponding sequence in $E$ such that $\beta^{(n)}\in\mathrm{Supp}\ (\alpha^{(n)})'$ for any $n$. We claim that, without loss of generality, the $\alpha^{(n)}$'s may be assumed to be pairwise distinct as in the other case. Indeed, since $(\beta^{(n)})_{n\in\mathbb{N}}$ is strictly decreasing and for any $\alpha$, $\mathrm{Supp}\ \alpha '$ is anti-well-ordered in $\Gamma$, we have $\{\beta^{(n)}\ |\ n\in\mathbb{N}\}\cap\mathrm{Supp}\ \alpha '$ is finite for any $\alpha$. Therefore the set $\{\alpha^{(n)}\ |\ n\in\mathbb{N}\}$ has to be infinite: it suffices to restrict to a subsequence of representatives of this set, which we continue to denote by $(\alpha^{(n)})_{n\in\mathbb{N}}$ below.
\sn
\textsl{From now on, we will not distinguish between the two preceding cases, writing that $\beta^{(n+1)}\succcurlyeq\beta^{(n)}$ for all $n$.} From (D1), we note that $\mathrm{Supp}\ \alpha '\subset\left(\alpha.\displaystyle\bigcup_{\phi\in\mathrm{supp}\ \alpha}\mathrm{Supp}\ \displaystyle\frac{\phi '}{\phi}\right)$  for any $\alpha$. Hence, for any $n$, we set
$\beta^{(n)}=\alpha^{(n)}\tau^{(n)}$ for some $\tau^{(n)}\in\mathrm{Supp}\ \displaystyle\frac{\phi_n'}{\phi_n}$ with $\phi_n\in\mathrm{supp}\ \alpha^{(n)}$. We now apply Lemma \ref{lemme:suite} to the sequence $S=(\alpha^{(n)})_{n\in\mathbb{N}}$ of pairwise distinct elements of $E$. Since  $E$ is anti-well-ordered in $\Gamma$, $S$ cannot have an infinite strictly increasing subsequence. So $S$ has a strictly decreasing subsequence which we continue to denote $(\alpha^{(n)})_{n\in\mathbb{N}}$ for convenience.
\sn
Since for any $k<l\in\mathbb{N}$,
$\beta^{(k)}=\alpha^{(k)}\tau^{(k)}\preccurlyeq\beta^{(l)}=\alpha^{(l)}\tau^{(l)}$,
we have : \begin{equation}\label{inequ1}\forall
k<l\in\mathbb{N},\ \
1\prec\displaystyle\frac{\alpha^{(k)}}{\alpha^{(l)}}\preccurlyeq
\displaystyle\frac{\tau^{(l)}}{\tau^{(k)}}
\end{equation}
The sequence $(\tau^{(n)})_{n\in\mathbb{N}}$ is therefore strictly increasing.\sn
Now consider a corresponding sequence $(\phi_{n})_{n\in\mathbb{N}}$ (for which $\tau^{(n)}\in\mathrm{Supp}\ \displaystyle\frac{\phi_{n}'}{\phi_{n}}$ and $\phi_n\in\mathrm{supp }\alpha^{(n)}$ for any $n$). As for the first case here above, we may assume without loss of generality that the $\phi_{n}$'s are pairwise distinct. 
\sn
We apply Lemma \ref{lemme:suite} to the sequence $\tilde{S}=(\phi_n)_{n\in\mathbb{N}}$. Suppose that it has an infinite decreasing subsequence, say $\hat{S}=(\phi_{n_i})_{i\in\mathbb{N}}$. This anti-well-ordered subset $\hat{S}\subset\Phi$ would be such that the corresponding subsequence $(\tau^{(n_i)})_{i\in\mathbb{N}}$ is increasing, contradicting $(SD1)$. So $\tilde{S}$ has an infinite increasing subsequence which we continue to denote $(\phi_n)_{n\in\mathbb{N}}$ for convenience.
\sn
We shall define by induction strictly increasing subsequences $(\phi_{n_i})_{i\in\mathbb{N}}$ of $(\phi_n)_{n\in\mathbb{N}}$ and $(\tau^{(n_i)})_{i\in\mathbb{N}}$ of $(\tau^{(n)})_{n\in\mathbb{N}}$ as in the statement of $(H2)$. Set $n_0=0$ and recall that for any $n$, $\phi_n\in\mathrm{supp}\  \alpha^{(n)}$. Suppose that we have subsequences $\phi_{n_0}\prec\phi_{n_1}\prec\cdots\prec\phi_{n_i}$ and $\tau^{(n_0)}\prec\tau^{(n_1)}\prec\cdots\prec\tau^{(n_i)}$ for some $i\geq 0$.
Since the sequence $(\phi_n)_{n\in\mathbb{N}}$ is increasing and
$\mathrm{supp}\ \alpha^{(n_i)}$ is anti-well-ordered in $\Phi$, there
exists a lowest index $n_{i+1}>n_i$ such that $\phi_{n_{i+1}}\notin\mathrm{supp}\
\alpha^{(n_i)}$. But $\phi_{n_{i+1}}\in\mathrm{supp}\ \alpha^{(n_{i+1})}$.
So $\phi_{n_{i+1}}\in\mathrm{supp}\ \displaystyle\frac{\alpha^{(n_i)}}{\alpha^{(n_{i+1})}}$ and $\mbox{ LF
}\left(\displaystyle\frac{\alpha^{(n_i)}}{\alpha^{(n_{i+1})}}\right) \succcurlyeq\phi_{n_{i+1}}$.
Moreover by (\ref{inequ1}) we have $\mbox{ LF
}\left(\displaystyle\frac{\tau^{(n_{i+1})}}{\tau^{(n_i)}}\right)\succcurlyeq \mbox{ LF
}\left(\displaystyle\frac{\alpha^{(n_i)}}{\alpha^{(n_{i+1})}}\right)$. So
$\mbox{ LF
}\left(\displaystyle\frac{\tau^{(n_{i+1})}}{\tau^{(n_i)}}\right)\succcurlyeq \phi_{n_{i+1}}$ as required.
\end{demonstration}

From Lemma \ref{lemme:H2} and Corollary \ref{coro:(H1')} we deduce a
more explicit sufficient (but not necessary) condition such that a
map $d\ :\ \Phi\rightarrow\mathds{K}\backslash\{0\}$ extends to a
series derivation on $\mathds{K}$:
\begin{cor}\label{coro:deriv'}
Consider a map $d\ :\ \Phi\rightarrow\mathds{K}\backslash\{0\}$. Then $d$ extends to a series derivation on $\mathds{K}$ if the following properties hold:
\begin{description}
\item[(H1')] $E_1:=\left\{\phi\in\Phi\ |\ \exists\ \psi\succ\phi,\
I_{\psi,\phi}\mathrm{\ is\ not\ a\ left\ shift}\right\}$ is well ordered
in $\Phi$.
\item[(H2'')] $E_2:= \left\{\psi\in\Phi\>|
\exists\ \phi\prec\psi,\> \exists \tau^{(\phi)}\in\mathrm{Supp}\>
\displaystyle\frac{\phi '}{\phi},\>
\exists \tau^{(\psi)}\in\textrm{Supp}\ \displaystyle\frac{\psi
'}{\psi}\> \mathrm{\ s.t.\ } \mbox{ LF
}\>(\displaystyle\frac{\tau^{(\phi)}}{\tau^{(\psi)}}\>)\succcurlyeq\psi
 \right\}$ is anti-well-ordered in $\Phi$.
\end{description}
\end{cor}
\begin{demonstration}
By Corollary \ref{coro:(H1')}, $d$ extends to a series derivation on
$\Gamma$. From Lemma \ref{lemme:H2}, (SD2) does not hold if and
only if there exist infinite increasing sequences
$(\phi_n)_{n\in\mathbb{N}}\subset\Phi$ and
$(\tau^{(n)})_{n\in\mathbb{N}}\subset\Gamma$ such that for any $n$,
$\tau^{(n)}\in\mathrm{Supp}\ \displaystyle\frac{\phi_n'}{\phi_n}$
and $\mbox{ LF
}\left(\displaystyle\frac{\tau^{(n+1)}}{\tau^{(n)}}\right)\succcurlyeq\phi_{n+1}$.
But from (H2''), for any increasing sequence
$S=(\phi_n)_{n\in\mathbb{N}}$, since $E_2\subset\Phi$ is
anti-well-ordered, $E_2\cap S$ is finite. So, for all but finitely
many $n$, $\mbox{ LF
}\left(\displaystyle\frac{\tau^{(n+1)}}{\tau^{(n)}}\right)\prec\phi_{n+1}$
for any $\tau^{(n)}\in\mathrm{Supp}\ \displaystyle\frac{\phi_n
'}{\phi_n}$ and any $\tau^{(n+1)}\in\mathrm{Supp}\
\displaystyle\frac{\phi_{n+1} '}{\phi_{n+1}}$. This contradicts
$(H2)$. \end{demonstration}

\begin{ex}
If we omit the assumption that the sequence
$(\tau^{(n)})_{n\in\mathbb{N}}$ is increasing in (H2) (or
(H2')), the condition is not anymore necessary, even if we
restrict to the case that the supports of $\displaystyle\frac{\phi
'}{\phi}$ are finite and uniformly bounded as in Corollary \ref{coro:(H1'')}. Indeed we have the following example. Given an infinite increasing sequence $(\phi_n)_{n\in\mathbb{N}}$, suppose that there exists $\psi\in\Phi$ such that $\psi\succ\phi_n$ for any
$n$. Then define
$\displaystyle\frac{\phi_0'}{\phi_0}=\tau^{(0)}_1+\tau^{(0)}_2=
1+\psi^{-1}$ and for any $n\in\mathds{N}^*$,
$\displaystyle\frac{\phi_n'}{\phi_n}=\tau^{(n)}_1+\tau^{(n)}_2= \phi_{n-1}+\psi^{-1}\phi_{n-1}$.

 We observe that any infinite increasing sequence of $\tau$'s contains either
infinitely many $\tau_1$'s, or infinitely $\tau_2$'s. Moreover
for any $k<l$, $\mbox{ LF
}\left(\displaystyle\frac{\tau^{(k)}_{1}}{\tau^{(l)}_{1}}\right)=
\mbox{ LF
}\left(\displaystyle\frac{\tau^{(k)}_{2}}{\tau^{(l)}_{2}}\right)= \phi_{l-1}\prec\phi_l$.
So  (SD2) holds, even if for any $n\in\mathbb{N}$, $\mbox{ LF
}\left(\displaystyle\frac{\tau^{(n+1)}_{2}}{\tau^{(n)}_{1}}\right)=
\mbox{ LF }\left(\displaystyle\frac{\tau^{(n+1)}_{1}}{\tau^{(n)}_{2}}\right)= \psi\succ\phi_{n+1}$.
\end{ex}

\section{Hardy type derivations.}\label{sect:ex}
\begin{defn}\label{defi:hardy_deriv} Let $(K,\preccurlyeq,\mathcal{C})$ be a field endowed with a dominance relation (cf. Definition \ref{defi:dom_rel}), which contains a sub-field $\mathcal{C}$ isomorphic to its
 residue field $K^{\preccurlyeq 1}/K^{\prec 1}$ (so
$K^{\preccurlyeq 1}=\mathcal{C}\oplus K^{\prec 1}$). A derivation
$d\ :\ K\rightarrow K\ $ is a {\bf Hardy type derivation} if :
\begin{description}
\item[(HD1)] the \textbf{sub-field of constants} of
$K$ is $\mathcal{C}$ : $\forall a\in K,\ a'=0\Leftrightarrow
a\in\mathcal{C}$.
\item[(HD2)] $d$ verifies \textbf{l'Hospital's rule} :
$\forall a,b\in K\backslash\{0\}$ with $a,b \nasymp 1$ we have\n
$a\preccurlyeq b\Leftrightarrow a'\preccurlyeq b'$.
\item[(HD3)] the logarithmic derivation is \textbf{compatible with the dominance relation} (in the sense of Hardy fields): \emph{$\forall  a,b\in K\ \ with\ |a|\succ|b|\succ 1,
we\ have\ \displaystyle\frac{a'}{a}\succcurlyeq
\displaystyle\frac{b'}{b}$. Moreover,
$\displaystyle\frac{a'}{a}\asymp\displaystyle\frac{b'}{b}$ if and
only if $a$ and $b$ are comparable.}
\end{description}
\end{defn}
Axioms (HD1) and (HD2) are exactly those which define a \textbf{differential
valuation} (\cite[Definition p. 303]{rosenlicht:diff_val}; see
Theorem 1 and Corollary 1 for the various versions of l'Hospital's
rule that hold in this context). Axiom (HD3) is the version for dominance relations of the Principle (*) in \cite[p. 992]{rosenlicht:val_gr_diff_val2}. This principle is itself a generalisation of properties obtained in \cite[Propositions 3 and 4]{rosenlicht:rank} and \cite[Principle (*) p.
314]{rosenlicht:diff_val} in the context of Hardy fields: recall that a \textbf{Hardy field} is, by definition, a field of germs at $\infty$ of real functions closed under differentiation \cite[Chap.V, App.]{bour:unevar}. E.g., the fields (of the corresponding germs) of real rational functions $\mathbb{R}(x)$, of real meromorphic functions at $+\infty$, of Logarithmic-Exponential functions \cite{hardy}\cite{hardy:LEfunct}. They are prime examples of differential valued field, the valuation being the natural one induced by the ordering of germs \cite{rosenlicht:hardy_fields}. 
\mn Below we prove the following criterion for a series derivation to be of Hardy type.
\begin{nota}\label{nota:theta}
Let $\phi\in\Phi$. Set $\theta^{(\phi)}\> := \> \mbox{ LM }(\phi '/\phi)$, i.e. \begin{center}
$\displaystyle\frac{\phi '}{\phi}=t_\phi\theta^{(\phi)}(1+\epsilon)$
\end{center}
where $t_\phi\in\mathbb{R}^*$ and $\epsilon\in \mathds{K}^{\prec 1}$.
\end{nota}
\begin{thm}\label{theo:hardy-deriv} A series
derivation $d$ on $\mathds{K}$ verifies \emph{(HD2)} and \emph{(HD3)} if and only if
the following condition holds:\begin{description}
\item[(H3')] $\forall\phi\prec\psi\in\Phi$,
$\theta^{(\phi)}\prec\theta^{(\psi)}$ and $\mbox{ LF
}\left(\displaystyle\frac{\theta^{(\phi)}}{\theta^{(\psi)}}\right)\prec\psi$.
\end{description}
\end{thm}
\begin{demonstration}
We suppose that (H3') holds. To prove l'Hospital's rule on
$\mathds{K}$, it suffices to prove it for the monomials. Let
$\alpha=\displaystyle\prod_{\phi\in\mathrm{supp}\
\alpha}\phi^{\alpha_\phi}$ and
$\beta=\displaystyle\prod_{\phi\in\mathrm{supp}\
\beta}\phi^{\beta_\phi}$ be arbitrary monomials. Then $\alpha
'=\alpha\displaystyle\sum_{\phi}\alpha_\phi\displaystyle\frac{\phi
'}{\phi}\asymp\alpha\theta^{(\phi_0)}$ and $\beta
'=\beta\displaystyle\sum_{\phi}\beta_\phi\displaystyle\frac{\phi
'}{\phi}\asymp\beta\theta^{(\phi_1)}$ where $\phi_0=\mbox{ LF
}(\alpha)$
and $\phi_1=\mbox{ LF }(\beta)$.

If $\phi_0=\phi_1$, then $\theta^{(\phi_0)}=\theta^{(\phi_1)}$.
So $\displaystyle\frac{\alpha '}{\beta '}\asymp \displaystyle\frac{\alpha}{\beta}$. If $\phi_0\neq\phi_1$, for instance $\phi_0\prec\phi_1$, then
$\mbox{ LF }\left(\displaystyle\frac{\alpha}{\beta}\right)=\phi_1$.
But $\displaystyle\frac{\alpha '}{\beta '}\asymp
\displaystyle\frac{\alpha}{\beta}
\displaystyle\frac{\theta^{(\phi_0)}}{\theta^{(\phi_1)}}$, and
$\mbox{ LF
}\left(\displaystyle\frac{\theta^{(\phi_0)}}{\theta^{(\phi_1)}}\right)\prec\phi_1$.
Applying the ultrametric inequality for LF, we obtain $\mbox{ LF }\left(\displaystyle\frac{\alpha '}{\beta
'}\right)=\phi_1$ and $\mbox{ LE }\left(\displaystyle\frac{\alpha
'}{\beta
'}\right)=\mbox{ LE }\left(\displaystyle\frac{\alpha}{\beta}\right)$. Thus $\displaystyle\frac{\alpha '}{\beta '}$ and $\displaystyle\frac{\alpha}{\beta}$ have same sign. \sn
 To prove the compatibility of the logarithmic derivation, take $a,b\in\mathds{K}\ \ with\
|a|\succ |b|\succ 1$ and denote $\alpha=\mbox{ LM }(a)$,
$\beta=\mbox{ LM }(b)$, $\phi_0=\mbox{ LF }(a)=\mbox{ LF }(\alpha)$
and $\phi_1=\mbox{ LF }(b)=\mbox{ LF }(\beta)$. So we have $\mbox{
LM
}\left(\displaystyle\frac{a'}{a}\right)=\mbox{ LM
}\left(\displaystyle\frac{\alpha
'}{\alpha}\right)=\mbox{ LM
}\left(\displaystyle\frac{\phi_0'}{\phi_0}\right)=\theta^{(\phi_0)}$
and similarly $\mbox{ LM
}\left(\displaystyle\frac{b'}{b}\right)=\theta^{(\phi_1)}$ (Lemma \ref{lemme:LF}). Since
$|a|\succ |b|\succ 1$, we have $\phi_0\succcurlyeq\phi_1$. So $\theta^{(\phi_0)}\succcurlyeq\theta^{(\phi_1)}$ by (H3').
Moreover, $a$ and $b$ are comparable if and only if $\phi_0=\phi_1$,
which means that $\theta^{(\phi_0)}=\theta^{(\phi_1)}$.
\mn Conversely, for $\phi, \psi \in \Phi$ with $\phi \prec
\psi$ we have $\displaystyle{\frac{\phi '}{\phi}}\prec
\displaystyle{\frac{\psi '}{\psi}}$, since the logarithmic derivation is assumed to be compatible with the dominance relation (recall that $\phi\succ 1$ for any $\phi\in\Phi$ by construction). Thus $\mbox{ LM
}(\displaystyle{\frac{\phi '}{\phi}})\prec \mbox{ LM
}(\displaystyle{\frac{\psi '}{\psi}})$, that is $\theta ^{(\phi)}
\prec \theta ^{(\psi)}$, and $1\prec
\displaystyle{\frac{\theta ^{(\psi)}} { \theta ^{(\phi)}}}$.
\sn 
Now, for any reals $r<0$ and $s\ne 0$ and any $\phi$, $\psi
\in \Phi$ with $\phi \prec \psi$, we have $\psi^r \prec \phi^s$.
Differentiating both sides and applying
l'Hospital's rule, we obtain $1 \prec \displaystyle\frac{\psi
'/\psi}{\phi '/\phi}\prec \phi^s \psi^{-r} \>.$ Now $\mbox{ LF
}(\phi ^s \psi^{-r}) = \psi$ and $\mbox{ LE }(\phi^s \psi^{-r}) = -r
>0$. Thus $\mbox{ LF }(\displaystyle{\frac{\theta ^{(\psi)}} {
\theta ^{(\phi)}}})\prec \psi \>.$
\end{demonstration}
\begin{cor}
 A series derivation $d$ on $\mathds{K}$ which verifies \emph{(HD2)} and \emph{(HD3)} is a Hardy type derivation.
\end{cor}
\begin{demonstration}
By construction the field of coefficients $\mathbb{R}$ is included in the field of constants (see (D0), (D2)). Conversely, consider a
non-constant series $a=\displaystyle\sum_{\alpha\in\mathrm{Supp}\
a}a_\alpha \alpha \in
\mathds{K}\backslash\{0\}$ such that $a'=0$. By (D1), we have $a'=\displaystyle\sum_{\alpha\in\mathrm{Supp}\ a}a_\alpha \alpha '$.
Set $\alpha^{(0)}=\max((\mathrm{Supp}\ a)\backslash
\{1\})$. By l'Hospital's rule, we have
$(\alpha^{(0)})'\succ\alpha '$ for any $\alpha\in((\mathrm{Supp}\
a)\backslash\{\alpha^{(0)}\})$.
Thus we would have $(\alpha^{(0)})'=0$. But, setting $\phi_0=\mbox{ LF }(\alpha^{(0)})$, by (D1)
and (H3') we obtain $(\alpha^{(0)})'\asymp
\alpha^{(0)}\theta^{(\phi_0)}$ which is non zero. Thus $(\alpha^{(0)})'$ cannot be zero, neither do $a'$ : this contradicts the initial assumption.
\end{demonstration}
\begin{remark}\label{rem:H-fields}
\sn \ \ In \cite{vdd:asch:liouv-cl-H-fields} is developed the notion of \textbf{H-field}, which generalises the one of Hardy field. Indeed, by definition,  an H-field is an ordered differential field endowed with a dominance relation $(K,d,\leq,\preccurlyeq)$ and with sub-field of constants $C$, such that the two following properties hold:
\begin{description}
    \item[(HF1)] if $f\succ 1$, then $\displaystyle\frac{f'}{f} >0$;
\item[(HF2)] if $f\preccurlyeq 1$, then $f-c\prec 1$ for some $c\in C$.
\end{description}
Therefore, in our context of generalised series endowed with a  Hardy type derivation, we note that (HD1) is equivalent to (HF2). Therefore:
\begin{center}
\emph{ $\mathds{K}$ is an H-field if and only if for any $\phi\in\Phi$, $\displaystyle\frac{\phi '}{\phi} >0$, i.e. $\mbox{LC }\left(\displaystyle\frac{\phi '}{\phi}\right)>0$. }
\end{center}
\noindent Indeed, for any series $a\succ 1$, denote $\mbox{ LM }(a)=\alpha$, $\mbox{ LF }(a)=\phi$ and  $\mbox{ LC }(a)=\alpha_0>0$. As noticed in the preceding proof, by (D1), (D2) and (H3'), we have: 
\begin{equation}\label{equ:deriv-log}
\displaystyle\frac{a'}{a}\sim\displaystyle\frac{\alpha '}{\alpha}\sim \alpha_0\displaystyle\frac{\phi '}{\phi}.
\end{equation}
So $\displaystyle\frac{a'}{a}$ has same sign as $\displaystyle\frac{\phi '}{\phi}$.
\end{remark}

\section{Examples.} \label{sect:exemple}

\subsection{The monomial case.}
\begin{defn}\label{defi:monom-deriv} A series derivation on $\mathds{K}$
is \textbf{monomial} if its restriction to the fundamental monomials has its image in the monomials: \begin{center}
$d:\Phi\rightarrow \mathbb{R}^*.\Gamma$
\end{center}
i.e. with Notation \ref{nota:theta}, we have \begin{center}
$\displaystyle\frac{\phi '}{\phi}:=t_\phi\theta^{(\phi)}$ for some $t_\phi \in  \mathbb{R}^*$.
\end{center}
\end{defn}

\begin{prop}\label{prop:monomial_case}
 A map $d:\Phi\rightarrow\mathbb{R}^*.\Gamma$ extends to a series derivation of Hardy type on $\mathds{K}$ if and only if the Hypothesis (H3') holds.
\end{prop}
\begin{demonstration}
Given a map $d:\Phi\rightarrow\mathbb{R}^*.\Gamma$, there exists a series derivation on $\mathds{K}$ (extending it) if and only if (H1'') with $N=1$ and (H2') hold (see Corollary \ref{coro:(H1'')} and Remark \ref{rem:H2'}). Then, it suffices to remark that (H3') is a particular case of (H1'') and (H2'), in which the only element in $\mathrm{Supp}\ \displaystyle\frac{\phi '}{\phi}$ is $\theta^{(\phi)}$. Now apply Theorem \ref{theo:hardy-deriv}.
\end{demonstration}

\begin{defn}
Given a totally ordered set $(\Phi,\preccurlyeq)$, we call a \textbf{left-shift endomorphism} of $\Phi$ any order preserving map $s : \Phi\rightarrow \Phi$ (i.e. $\phi_1\prec\phi_2\Leftrightarrow s(\phi_1)\prec s(\phi_2)$) such that $s(\phi)\prec\phi$ for any $\phi\in \Phi$. Note that this implies that $\Phi$ has no least element. For any $n\in\mathbb{N}$, we denote by $s^n$ the $n^{\mathrm{th}}$ iterate of $s$.
\end{defn}

\textbf{1.} Let $(\Phi,\preccurlyeq)$ be a totally ordered set that we suppose endowed with a left-shift endomorphism $s\ :\ \Phi\rightarrow \Phi$. Set $\theta^{(\phi)}:=s(\phi) \mathrm{\ \ for\ any\ } \phi\in\Phi$. We claim that for any choice of $t_\phi\in\mathbb{R}^*$, the corresponding map \[\begin{array}{llcl}
d:&\Phi&\rightarrow&\mathbb{R}^*.\Gamma\\
&\phi&\mapsto& t_\phi.s(\phi)\phi
\end{array}\]
  extends to a series derivation of Hardy type. Indeed, by Proposition \ref{prop:monomial_case}, it suffices to show that Hypothesis (H3') holds. Indeed, for any $\phi_1\neq \phi_2$, we have
$\phi_1\prec\phi_2\Leftrightarrow\theta^{(\phi_1)} =s(\phi_1)\prec
\theta^{(\phi_2)}=s(\phi_2)$. Moreover, $\mbox{ LF
}\left(\displaystyle\frac{\theta^{(\phi_1)}}{\theta^{(\phi_2)}}\right)
= \mbox{ LF
}\left(\displaystyle\frac{s(\phi_1)}{s(\phi_2)}\right)=
s(\phi_2)\prec\phi_2$. That is, we have (H3').

Note that we could have set
$\theta^{(\phi)}:=s(\phi)^{\alpha_\phi}$ for some $\alpha_\phi>0$ (see Proposition \ref{prop:ex_gene}). 

Note also that the preceding example extends to the case where $\Phi$ has a least element $\phi_m$ and $\Phi\setminus\{\phi_m\}$ carries a left-shift endomorphism, just by setting $\theta^{(\phi_m)}:=1$ and $\theta^{(\phi)}:=s(\phi)$ for any $\phi\succ \phi_m$.
\mn

\textbf{2.} We generalise the preceding example. For any $\phi\in\Phi$, fix $N_\phi\in\mathbb{N}\cup\{+\infty\}$. One can set $\theta^{(\phi)}:=\displaystyle\prod_{n=1}^{N_\phi}s^n(\phi)$ for any $\phi\in \Phi$. As above, Hypothesis (H3') holds.\\
Note that we can also set $\theta^{(\phi)}:=\displaystyle\prod_{n=1}^{N_\phi}s^n(\phi)^{\alpha_{\phi,n}}$ with $\alpha_{\phi,n}\in\mathbb{R}$ for all $n\in\mathbb{N}$ and  $\alpha_{\phi,1}>0$ (see Proposition \ref{prop:ex_gene}).
%
%
\mn

\textbf{3.} Assume now that $\Phi$ is isomorphic to a subset of
$\mathbb{R}$ with least element $\phi_m$, writing $f$ this
isomorphism, we can set for any $\phi\in \Phi$,
$\theta^{(\phi)}:=\phi_m^{f(\phi)+\beta}$ where
$\beta$ is some fixed real.
\sn 
As an illustration, take  the following chain of infinitely increasing real germs at infinity (applying the usual comparison relations of germs) 
$\Phi=\{\phi_\alpha=\exp(x^\alpha)\ ;\ \alpha>0\}\cup\{\phi_0=x\}$
which is isomorphic to $\mathbb{R}_+$. With the usual derivation of (germs at infinity of) real functions, we have
$\phi_\alpha ' =\alpha
x^{\alpha-1}\exp(x^\alpha)=\alpha \phi_0^{\alpha-1}\phi_\alpha$ and
$\phi_0'=1$. Thus, $\theta^{(\phi_\alpha)}=\phi_0^{\alpha-1}$ and
$t_{\phi_\alpha}=\alpha$. 
\mn

\textbf{4.} Assume that $\Phi$ is anti-well-ordered, with greatest element $\phi_M$. Consider some fixed \textit{negative} reals $\alpha_\psi$ for $\psi\in\Phi$. We can set $\theta^{(\phi)}:=\displaystyle\prod_{\psi\prec\phi} \psi^{\alpha_{\psi,\phi}}
\displaystyle\prod_{\phi\succcurlyeq \psi\succcurlyeq \phi_M}\psi^{\alpha_\psi}$ with arbitrary $\alpha_{\psi,\phi}\in\mathbb{R}$, provided that   $\alpha_{\psi_\phi,\phi}>\alpha_{\psi_\phi}$ where $\psi_\phi$ denotes the predecessor of $\phi$ in $\Phi$.
\mn
%
%
%
As an illustration of examples 2 and 4, take now $\Phi:=\{\exp^n(x)\ ;\ n\in\mathbb{Z}\}$
where $\exp^n$ denotes for positive $n$, the $n$'th iteration of the real exponential function, for negative $n$,
the $|n|$'s iteration of the logarithmic function, and for $n=0$ the identical map. Note that $\mathds{K}$ contains naturally the field of rational fractions $\mathbb{R}(\exp^n(x),\ n\in\mathbb{Z})$ which is a Hardy field (see the commentaries after Definition \ref{defi:hardy_deriv}). We have: \begin{center}
$\left\{\begin{array}{lclcll}
\displaystyle\frac{(\exp^n(x))'}{\exp^n(x)} &=&\theta^{(n)}(x)&=&\displaystyle\prod_{k=1}^{n-1} \exp^k(x)\ &\textrm{if}\ n\geq 2\\
\displaystyle\frac{(\exp(x))'}{\exp(x)}&=&\theta^{(1)}(x) &=&1\ &\\
\displaystyle\frac{(\exp^n(x))'}{\exp^n(x)}&=&\theta^{(n)}(x) &=&\displaystyle\prod_{k=0}^{n} \displaystyle\frac{1}{\exp^k(x)}\ &\textrm{if}\ n\leq 0
\end{array}\right.$
\end{center}
So for any integers $m<n$, we have  $\theta^{(m)}\prec \theta^{(n)}$ and $\mbox{LF }\left(\displaystyle\frac{\theta^{(n)}}{\theta^{(m)}}\right) =\exp^{n-1}(x)\prec \exp^n(x)$: (H3') holds. By Proposition \ref{prop:monomial_case}, the usual derivation of germs $\exp^n(x)\mapsto(\exp^n(x))'$ extends to a series derivation of Hardy type on $\mathds{K}$. Moreover, since the leading coefficients of $\displaystyle\frac{(\exp_n(x))'}{\exp_n(x)}$ is always 1 which is positive, $\mathds{K}$ endowed with such a derivation is an H-field (see Remark \ref{rem:H-fields}).

\subsection{A general example.}
To motivate the introduction of the non monomial case, consider the  Hardy field\\ $\mathbb{R}(x,\exp(x),\exp(x^2),\exp(\exp(x^2+x)))$ (for $x$ near $\infty$). Then, denoting $\phi=\exp(\exp(x^2+x))$, we have
\begin{center}
$\displaystyle\frac{\phi '}{\phi}=2x\exp(x^2)\exp(x)+\exp(x^2)\exp(x)$
\end{center}
which is not a monomial.

We proceed by generalizing the preceding examples. We suppose that $\Phi$ carries a left-shift endomorphism $s:\Phi\rightarrow\Phi$. We shall define a family of derivations on $\mathds{K}$. This family is defined using the following other field of generalised series.\sn
We consider an ordered set of fundamental monomials $(\Lambda=\{\lambda_n\ ;\ n\in\mathbb{N}\},\preccurlyeq)$ isomorphic to $(\mathbb{N},\leq)$, the corresponding Hahn group of monomials $\mathrm{H}(\Lambda)$ and field of generalised series $\mathbb{L}:=\mathbb{R}((\mathrm{H}(\Lambda)))$ as in Section \ref{sect:defi}. We recall that $\mathbb{L}^{\succ 1}$ denotes the subring of purely infinite series, which is an additive complement group of the valuation ring in $\mathbb{L}$.

\begin{prop}\label{prop:ex_gene}
For any purely infinite series $l=\displaystyle\sum_{\delta\in S}l_\delta \displaystyle\prod_{n\in\mathbb{N}}\lambda_n^{\delta_n}$ (where $S$ denotes the support of $l$ which is in $\mathrm{H}(\Lambda)^{\succ 1}$), for any $\gamma\in\Gamma$, the map:\[d_{l,\gamma}\ :\ \phi\mapsto\gamma\phi\displaystyle\sum_{\delta\in S}l_\delta \displaystyle\prod_{n\in\mathbb{N}}[s^{n+1}(\phi)]^{\delta_n}\]
 is well-defined with values in $\mathds{K}$ (where $s^{n+1}$ denotes the $(n+1)^{\mathrm{th}}$ iterate of $s$), and it extends to a series derivation of Hardy type on $\mathds{K}$.
\end{prop}
\begin{demonstration}
We prove that conditions (H1') and (H2'') of Corollary \ref{coro:deriv'} and (H3') of Theorem \ref{theo:hardy-deriv} hold. Note that for any $\phi\in\Phi$, we have: 
\[\displaystyle\frac{\phi '}{\phi}=\displaystyle\sum_{\delta\in S}l_\delta.\gamma\displaystyle\prod_{n\in\mathbb{N}} [s^{n+1}(\phi)]^{\delta_n}.\]
For any $\phi\succ\psi$ in $\Phi$, the ordered sets $\mathrm{Supp}\
\displaystyle\frac{\phi '}{\phi}$ and $\mathrm{Supp}\
\displaystyle\frac{\psi '}{\psi}$ are isomorphic by construction.
Moreover, consider some monomial $\tau^{(\phi)}\in\mathrm{Supp}\
\displaystyle\frac{\phi '}{\phi}$, say
$\tau^{(\phi)}=\gamma\displaystyle\prod_{n\in\mathbb{N}} [s^{n+1}(\phi)]^{\delta_n}$
for some real $\delta_n$'s, $n\in\mathbb{N}$. Then we have
$I_{\phi,\psi}(\tau^{(\phi)})=\tau^{(\psi)}$ where
$\tau^{(\psi)}=\gamma\displaystyle\prod_{n\in\mathbb{N}} [s^{n+1}(\psi)]^{\delta_n}$.
Moreover, $\displaystyle\frac{\tau^{(\phi)}}{\tau^{(\psi)}}=
\displaystyle\prod_{n\in\mathbb{N}}\left(s^{n+1}(\phi)\right)^{\tau_n}
\left(s^{n+1}(\psi)\right)^{-\tau_n}$ with for all $n$,
$s^{n+1}(\phi)\succ s^{n+1}(\psi)$ (since $s$ is an endomorphism). Thus
$\mbox{ LF
}\left(\displaystyle\frac{\tau^{(\phi)}}{\tau^{(\psi)}}\right)=s^{n_0+1}(\phi)$
for some $n_0\in\mathbb{N}$. Moreover $\mbox{ LE
}\left(\displaystyle\frac{\tau^{(\phi)}}{\tau^{(\psi)}}\right)=\delta_{n_0}$
which is positive (since $d\in\mathbb{L}^{\succ 1}$). Hence we
obtain that $\displaystyle\frac{\tau^{(\phi)}}{\tau^{(\psi)}}\succ 1$, which means that $I_{\phi,\psi}$ is a decreasing automorphism. The Condition (H1') holds (the set $E_1$ is empty).

Consider now any $\tau^{(\phi)}=\gamma\displaystyle\prod_{n\in\mathbb{N}} [s^{n+1}(\phi)]^{\delta_n}\in\mathrm{Supp}\
\displaystyle\frac{\phi '}{\phi} $  and $\tau^{(\psi)}=\gamma\displaystyle\prod_{n\in\mathbb{N}} [s^{n+1}(\psi)]^{\delta'_n} \in \mathrm{Supp}\
\displaystyle\frac{\psi '}{\psi}$. Then $\mbox{ LF }\left(\displaystyle\frac{\tau^{(\phi)}}{\tau^{(\psi)}}\right)$ is equal to some $s^{n_0+1}(\phi)$ or $s^{n_0+1}(\psi)$ which is always less than $\phi$ since $s$ is a decreasing emdomorphism of $\Phi$. The Condition (H2'') holds (the set $E_2$ is empty).

Finally, note that the same properties hold in particular for the leading monomials $\theta^{(\phi)}$ and  $\theta^{(\psi)}$ of $\displaystyle\frac{\phi '}{\phi}$ and $\displaystyle\frac{\psi '}{\psi}$. The condition (H3') holds.
\end{demonstration}
\section{Asymptotic integration and integration}

\begin{defn} Let $(K,d,\preccurlyeq)$ be a differential field endowed with a dominance relation $\preccurlyeq$, and let $a$ be one of its elements.
 We say that $a$ admits an \textbf{asymptotic integral} $b$ if there exists $b\in K\setminus\{0\}$ such that $b'-a\prec a$. We say that $a$ admits an \textbf{integral} $b$ if there exists $b\in K\setminus\{0\}$ such that $b'=a$.
\end{defn}
The following main result about asymptotic integration in fields endowed with a Hardy type derivation is an adaptation of \cite[Proposition 2 and Theorem 1]{rosenlicht:rank}.
\begin{thm}[Rosenlicht] Let $(K,\preccurlyeq,\mathcal{C},d)$ be a field endowed with a Hardy type derivation
$d$.
 Let $a\in K\backslash\{0\}$, then $a$ admits an asymptotic
integral if and only if \begin{center}
 $a\nasymp\ g.l.b._\preccurlyeq\left\{\displaystyle\frac{b'}{b};\ b\in K\backslash\{0\},\ b\nasymp 1\right\}$
\end{center}
Moreover, for any such $a$, there exists $u_0\in
K\backslash\{0\}$ with $u_0\nasymp\ 1$ such that for any $u\in
K\backslash\{0\}$ such that $|u_0|\succcurlyeq |u|\succ 1$,
then
\begin{center}
$\left(a.\displaystyle\frac{au/u'}{(au/u')'}\right)'\sim a$
\end{center}
\end{thm}
\begin{demonstration} Our statement is a straightforward combination of Proposition 2 and Theorem 1 in \cite{rosenlicht:rank}.
It suffices to observe that the corresponding proofs in \cite{rosenlicht:rank} only rely on the fact that the canonical valuation of a
Hardy field is a differential valuation and that the logarithmic derivation is compatible with the dominance relation \cite[Proposition 3]{rosenlicht:rank}.
\end{demonstration}

In \cite[Lemma 1]{rosenlicht:rank}, Rosenlicht provides a method to compute $u_0$:
\begin{itemize}
    \item  since $a\nasymp\
g.l.b._\preccurlyeq\left\{\displaystyle\frac{b'}{b};\ b\in
K\backslash\{0\},\ b\nasymp 1\right\}$, we assume w.l.o.g. that \\
$a\succ g.l.b._\preccurlyeq\left\{\displaystyle\frac{b'}{b};\ b\in
K\backslash\{0\},\ b\nasymp 1\right\}$ (if not, take $a^{-1}$ instead of $a$);
\item take $u_1\succ 1$ such that $a\succ\displaystyle\frac{u_1'}{u_1}$;
\item  take any $u_0$ such that $u_0^{\pm 1}\asymp\min\left\{u_1,\ \displaystyle\frac{a}{u'_1/u_1}\right\}$.
\end{itemize}
Note that $u_0$ verifies $1\succ u_0^{\pm 1}\succcurlyeq
\left(\displaystyle\frac{a}{u_0'/u_0}\right)^{\pm 1}$. So, in our context, $\mbox{ LF }(u_0)\preccurlyeq \mbox{ LF
}\left(\displaystyle\frac{a}{u_0'/u_0}\right)$.

Our contribution here is to deduce explicit formulas for asymptotic integrals for our field of generalised series $\mathds{K}=\mathbb{R}((\Gamma))$
endowed with a Hardy type derivation. Note that this is equivalent (by l'Hospital's rule) to provide formulas for asymptotic integrals of monomials.
By (\ref{equ:deriv-log}) in Remark \ref{rem:H-fields}, note that we have:
\begin{center}
$g.l.b._\preccurlyeq\left\{\displaystyle\frac{b'}{b};\
b\in\mathds{K}\backslash\{0\},\ b\nasymp\ 1\right\} =
g.l.b._\preccurlyeq\left\{\displaystyle\frac{\phi '}{\phi};\
\phi\in\Phi\right\}= g.l.b._\preccurlyeq\left\{\theta^{(\phi)};\
\phi\in\Phi\right\}$.
\end{center}

Recall that for any monomial $\alpha\in\Gamma$, $\alpha\nasymp 1$, and any $\psi\in\mbox{supp }\alpha$, $\alpha_\psi$ denotes the exponent of $\psi$ in $\alpha$.

\begin{cor}\label{coro:asymp-int-monomial} Let $\alpha\in\Gamma$ be some monomial such that
$\alpha\nasymp\ g.l.b._\preccurlyeq\left\{\theta^{(\phi)};\
\phi\in\Phi\right\}$. If $\alpha\nasymp 1$, set $\phi_0:=\mbox{ LF
}(\alpha)$, so $\mbox{ LT }\left(\displaystyle\frac{\alpha
'}{\alpha}\right)
=\alpha_{\phi_0} t_{\phi_0}\theta^{(\phi_0)}$ (Remark \ref{rem:H-fields} (\ref{equ:deriv-log})). Then we have:
\begin{itemize}
    \item if $\mbox{ LF }(\theta^{(\phi_0)})\preccurlyeq\phi_0\asymp
\mbox{ LF
}\left(\displaystyle\frac{\alpha}{\theta^{(\phi_0)}}\right)$,
 then $\left(\displaystyle\frac{1}{t_{\phi_0}(\alpha_{\phi_0}-\theta^{(\phi_0)}_{\phi_0})}
 .\displaystyle\frac{\alpha}{\theta^{(\phi_0)}}\right)'\sim\alpha$;
\item if $\mbox{ LF }(\theta^{(\phi_0)})=\phi_1\succ\phi_0$,
then
$\left(\displaystyle\frac{1}{-t_{\phi_1}\theta^{(\phi_0)}_{\phi_1}}
.\displaystyle\frac{\alpha}{\theta^{(\phi_1)}}\right)'\sim\alpha$
(note that $\theta^{(\phi_0)}_{\phi_1}=\theta^{(\phi_1)}_{\phi_1}$);
\item if $\mbox{ LF }(\theta^{(\phi_0)})\asymp\phi_0\succ \mbox{
LF }\left(\displaystyle\frac{\alpha}{\theta^{(\phi_0)}}\right)$ or if $\alpha=1$,
then
$\left(\displaystyle\frac{1}{t_{\phi_1}(\alpha_{\phi_1}-\theta^{(\phi_1)}_{\phi_1})}
.\displaystyle\frac{\alpha}{\theta^{(\phi_1)}}\right)'\sim\alpha$
where $\phi_1$ is the element of $\Phi$ such that $\mbox{ LF
}\left(\displaystyle\frac{\alpha}{\theta^{(\phi_1)}}\right)=\phi_1\prec\phi_0$.
\end{itemize} 
\end{cor}
\begin{demonstration}
For the first case, it suffices to observe that $\mbox{ LF
}\left(\displaystyle\frac{\alpha}{\theta^{(\phi_0)}}\right)=\phi_0$
with exponent $\alpha_{\phi_0}-\theta^{(\phi_0)}_{\phi_0}$. So we have:
\begin{center}
$\left(\displaystyle\frac{\alpha}{\theta^{(\phi_0)}}\right)'\sim
\displaystyle\frac{\alpha}{\theta^{(\phi_0)}}
(\alpha_{\phi_0}-\theta^{(\phi_0)}_{\phi_0})\displaystyle\frac{\phi_0'}{\phi_0}\sim
\displaystyle\frac{\alpha}{\theta^{(\phi_0)}}
(\alpha_{\phi_0}-\theta^{(\phi_0)}_{\phi_0})t_{\phi_0}\theta^{(\phi_0)}= 
(\alpha_{\phi_0}-\theta^{(\phi_0)}_{\phi_0})t_{\phi_0}.\alpha$.
\end{center}

For the second case, since $\mbox{ LF
}(\theta^{(\phi_1)})=\phi_1\succ\phi_0$, we deduce from (H3) that
$\mbox{ LF }(\theta^{(\phi_0)})=\phi_1$ with the same exponent
$\theta^{(\phi_0)}_{\phi_1}$. So $\mbox{ LF
}\left(\displaystyle\frac{\alpha}{\theta^{(\phi_1)}}\right)=\phi_1$
with exponent $-\theta^{(\phi_0)}_{\phi_1}$, and then:\begin{center}
$\left(\displaystyle\frac{\alpha}{\theta^{(\phi_1)}}\right)'\sim
\displaystyle\frac{\alpha}{\theta^{(\phi_1)}}(-\theta^{(\phi_0)}_{\phi_1})
\displaystyle\frac{\phi_1'}{\phi_1}\sim \displaystyle\frac{\alpha}
{\theta^{(\phi_1)}}(-\theta^{(\phi_0)}_{\phi_1}t_{\phi_1})\theta^{(\phi_1)}= -\theta^{(\phi_0)}_{\phi_1}t_{\phi_1}.\alpha$.
\end{center}

For the third case, first we show that there exists $\phi_1$ as in the statement of the corollary. We define $u_0$ corresponding to $\alpha$ as in the preceding theorem and we denote $\hat{\phi}_0=\mbox{ LF }(u_0)$ and $\phi_1=\mbox{ LF
}\left(\displaystyle\frac{\alpha}{u_0'/u_0}\right)$. So we have
$\mbox{ LM
}\left(\displaystyle\frac{u_0'}{u_0}\right)=\beta_0\theta^{(\hat{\phi}_0)}$. Moreover by Rosenlicht's computation of $u_0$, we note that $\hat{\phi}_0\preccurlyeq\phi_1$. Thus we obtain by  (H3) that $\mbox{ LF
}\left(\displaystyle\frac{\theta^{(\hat{\phi}_0)}}{\theta^{(\phi_1)}}\right)\prec\phi_1$.
 and as desired:
\begin{center}
$\phi_1=\mbox{ LF }\left(\displaystyle\frac{\alpha}{\theta^{(\hat{\phi}_0)}}\right)
=\mbox{ LF
}\left(\displaystyle\frac{\alpha}{\theta^{(\hat{\phi}_0)}}.
\displaystyle\frac{\theta^{(\hat{\phi}_0)}}{\theta^{(\phi_1)}}\right)
=\mbox{ LF }\left(\displaystyle\frac{\alpha}{\theta^{(\phi_1)}}\right)$
\end{center}
Now we compute:
\begin{center}
$\begin{array}{lcl}
\left(\displaystyle\frac{\alpha}{\theta^{(\phi_1)}}\right)'&\sim&
\displaystyle\frac{\alpha}{\theta^{(\phi_1)}} (\alpha_{\phi_1}-\theta_{\phi_1}^{(\phi_1)}) \displaystyle\frac{\phi_1'}{\phi_1}\\
&\sim&(\alpha_{\phi_1}-\theta_{\phi_1}^{(\phi_1)}) \displaystyle\frac{\alpha}{\theta^{(\phi_1)}}t_{\phi_1}\theta^{(\phi_1)}\\
&=&t_{\phi_1} (\alpha_{\phi_1}-\theta_{\phi_1}^{(\phi_1)}). \alpha
\end{array}$.\end{center}
\end{demonstration}
Concerning integration, we apply to our context \cite[Theorem 55]{kuhl:Hensel-lemma} (recall that fields of generalised series are pseudo-complete  (see
e.g. \cite[Theorem 4, p. 309]{kap}).
\begin{cor} \label{int}
Assume that $\mathds{K}$ is endowed with a series derivation of
Hardy type  $d$. Set
$\tilde{\theta}=g.l.b._\preccurlyeq\left\{\theta^{(\phi)}\ ;\
\phi\in\Phi\right\}$ (if it exists). Then any element
$a\in\mathds{K}$ with $a\prec\tilde{\theta}$ admits an integral in
$\mathds{K}$. Moreover $\mathds{K}$ is closed under integration if
and only if $\tilde{\theta}\notin\Gamma$.
\end{cor}
\begin{demonstration}
As was already noticed before the Corollary \ref{coro:asymp-int-monomial}:
\begin{center}
$g.l.b._\preccurlyeq\left\{\displaystyle\frac{b'}{b};\
b\in\mathds{K}\backslash\{0\},\ b\nasymp\ 1\right\}=
\tilde{\theta}$.
\end{center}
Given $a\in\mathds{K}$ with
$a\prec\tilde{\theta}$, there exists a monomial
$\gamma\in\Gamma$ which is an asymptotic integral of $a$. That is,
 $\gamma '\asymp a$. Since $d$ verifies l'Hospital's rule, it
implies that for any $\tilde{\gamma}\in\mathrm{Supp}\ \gamma '$,
$\tilde{\gamma}\prec\tilde{\theta}$. So it admits itself an
asymptotic integral. The result now follows from  \cite[Theorem 55]{kuhl:Hensel-lemma}.
\end{demonstration}
Examples 2 and 3 in the case where $\Phi$ has no least element and the one of Proposition \ref{prop:ex_gene} are closed under integration.


%
%
%
\bibliographystyle{jloganal}
\def\cprime{$'$}

\adresse
\end{document}